%% file: main.tex
\pdftrailerid{}
\documentclass[journal,twoside,web]{ieeecolor}
\usepackage{xcolor}
\usepackage{generic}
\usepackage{cite}

\usepackage{amssymb}
\usepackage{amsmath}
\usepackage{mathtools}
\usepackage[hidelinks]{hyperref}

\usepackage{amsthm}
\usepackage{graphicx}
\usepackage{wrapfig}
\usepackage{capt-of}
\graphicspath{{figures/}}
\usepackage{cuted}
\usepackage{booktabs}

\usepackage{tikz}
\usepackage{pgfplots}
\pgfplotsset{compat=1.18}
\usepgfplotslibrary{groupplots}
\usepackage{mdframed}
\usepackage{etoolbox}
\makeatletter
\@ifundefined{@IEEEfigurecaptionsepspace}
  {\def\@IEEEfigurecaptionsepspace{\vskip\abovecaptionskip\relax}}{}
\@ifundefined{@IEEEtablecaptionsepspace}
  {\def\@IEEEtablecaptionsepspace{\vskip\abovecaptionskip\relax}}{}
\@ifundefined{ifcenterfigcaptions}{\newif\ifcenterfigcaptions}{}
\def\@IEEEtablestring{table}
\definecolor{captionlabelcolor}{cmyk}{1,0.53,0,0.47}
\long\def\@makecaption#1#2{%
  \ifx\@captype\@IEEEtablestring
    \begin{center}{\footnotesize\color{captionlabelcolor}#1}\\%
    {\footnotesize\scshape #2}\end{center}%
    \@IEEEtablecaptionsepspace
  \else
    \@IEEEfigurecaptionsepspace
    \setbox\@tempboxa\hbox{\def\color##1##{\@gobble}\footnotesize #1.~~ #2}%
    \ifdim\wd\@tempboxa>\hsize
      \parbox[t]{\hsize}{\footnotesize\noindent{\color{captionlabelcolor}#1.}~~ #2}%
    \else
      \ifcenterfigcaptions
        \hbox to\hsize{\footnotesize\hfil{\color{captionlabelcolor}#1.}~~ #2\hfil}%
      \else
        \hbox to\hsize{\footnotesize{\color{captionlabelcolor}#1.}~~ #2\hfil}%
      \fi
    \fi
  \fi}
\makeatother
\AtEndEnvironment{figure}{\normalcolor}
\AtEndEnvironment{figure*}{\normalcolor}
\AtEndEnvironment{table}{\normalcolor}
\AtEndEnvironment{theorem}{\normalcolor}
\makeatletter
\renewcommand{\fps@figure}{!t}
\renewcommand{\fps@table}{!t}
\makeatother
\definecolor{subsectioncolor}{cmyk}{1,0.53,0,0.47}
\definecolor{jinfocolor}{cmyk}{1,0.53,0,0.47}
\definecolor{boxbg}{rgb}{0.95, 0.95, 1.0}
\definecolor{boxborder}{rgb}{0.2, 0.2, 0.7}
\mdfdefinestyle{theorembox}{%
    topline=false, bottomline=false, rightline=false,
    linewidth=2pt, linecolor=boxborder,
    backgroundcolor=boxbg,
    innerleftmargin=5pt, innerrightmargin=5pt,
    innertopmargin=5pt, innerbottommargin=5pt
}

\newcommand{\reals}{\mathbb{R}}
\newcommand{\naturals}{\mathbb{N}}
\newcommand{\norm}[1]{\Vert #1 \Vert}
\newcommand{\abs}[1]{| #1 |}
\newcommand{\figscale}{1.00}

\def\logoname{blank_logo}

\newtheorem{theorem}{Theorem}
\newtheorem{lemma}[theorem]{Lemma}
\newtheorem{proposition}[theorem]{Proposition}
\newtheorem{corollary}[theorem]{Corollary}
\newtheorem{definition}{Definition}
\newtheorem{assumption}{Assumption}
\newtheorem{remark}{Remark}
\newtheorem{example}{Example}
\surroundwithmdframed[style=theorembox]{theorem}

\def\BibTeX{{\rm B\kern-.05em{\sc i\kern-.025em b}\kern-.08em
    T\kern-.1667em\lower.7ex\hbox{E}\kern-.125emX}}

\usepackage{eso-pic}
\input{posting_notice}

\begin{document}

\title{Strategically Robust Multi-Agent Control}

\author{Boris Velasevic, Nicolas Lanzetti, Eric Mazumdar
\thanks{%
NL is supported by the PIMCO Fellows Program, by the National Science Foundation under Grant~\#CCF-2326609, and by the Resnick Sustainability Institute.
EM is supported by National Science Foundation Award 2240110.
}%
\thanks{All authors are with the Department of Computing and Mathematical Sciences, California Institute of Technology, Pasadena, USA, (e-mail: \{bvelasev, lnicolas, mazumdar\}@caltech.edu).}}

\maketitle

\begin{abstract}
A key source of uncertainty in multi-agent systems is the behavior of other agents. Game-theoretic control approaches ignore this uncertainty, while robust control approaches treat others as fully adversarial.
In this paper, we introduce strategically robust multi-agent control to populate the spectrum between these two paradigms.
Specifically, in dynamic games, we propose that players are \emph{strategically} robust and seek protection against misspecified control policies of the other players.
This formulation induces a novel equilibrium concept, \emph{strategically robust dynamic equilibrium}, that interpolates between game-theoretic control and robust control. We study this concept across canonical open- and closed-loop settings and show that accounting for strategic robustness also has computational benefits over classic approaches in game-theoretic control. Numerical experiments on vehicle platooning, inventory management, and network control illustrate resilience to misspecification, and also exhibit a ``free-lunch'' phenomenon in which strategic robustness can improve performance or social welfare.

\end{abstract}

\begin{IEEEkeywords}
Multi-agent control, Dynamic games, Game theory, Robust control, Optimal control.
\end{IEEEkeywords}

\section{Introduction}\label{sec:intro}

\IEEEPARstart{I}n multi-agent systems like autonomous driving~\cite{abdel2024matched}, electricity markets~\cite{ruiz2009pool}, and real-world robotics tasks~\cite{fisac2018probabilistically}, one of the largest sources of uncertainty is the behavior of others. This \emph{strategic} uncertainty can arise from various sources such as misspecified models of others' intent, irrational behavior, or limited computation. Unlike standard exogenous disturbances commonly studied in control, strategic uncertainty is endogenous to the problem: by changing their behavior, an agent affects the uncertainty they are trying to address. 

Two common approaches in multi-agent control are to ignore strategic uncertainty and compute game-theoretic controllers under the implicit assumption that other agents’ objectives are well-modeled, or to account for strategic uncertainty through adversarial or worst-case modeling, i.e., robust control.

\begin{figure}[!t]
    \centering
    \input{figures/meta_figure}
    \vspace{-5pt}
    \caption{Strategically robust multi-agent control at a glance. Each player chooses how robust to be. Without robustness, players trust each other’s policies, recovering game-theoretic control; as robustness grows, each player hedges against ever larger deviations of the others, approaching robust control. Strategically robust equilibria fill the space in between and, once players are sufficiently robust (green), they are unique and efficiently computable.}
    \label{fig:meta}
    \vspace{-15pt}
\end{figure}
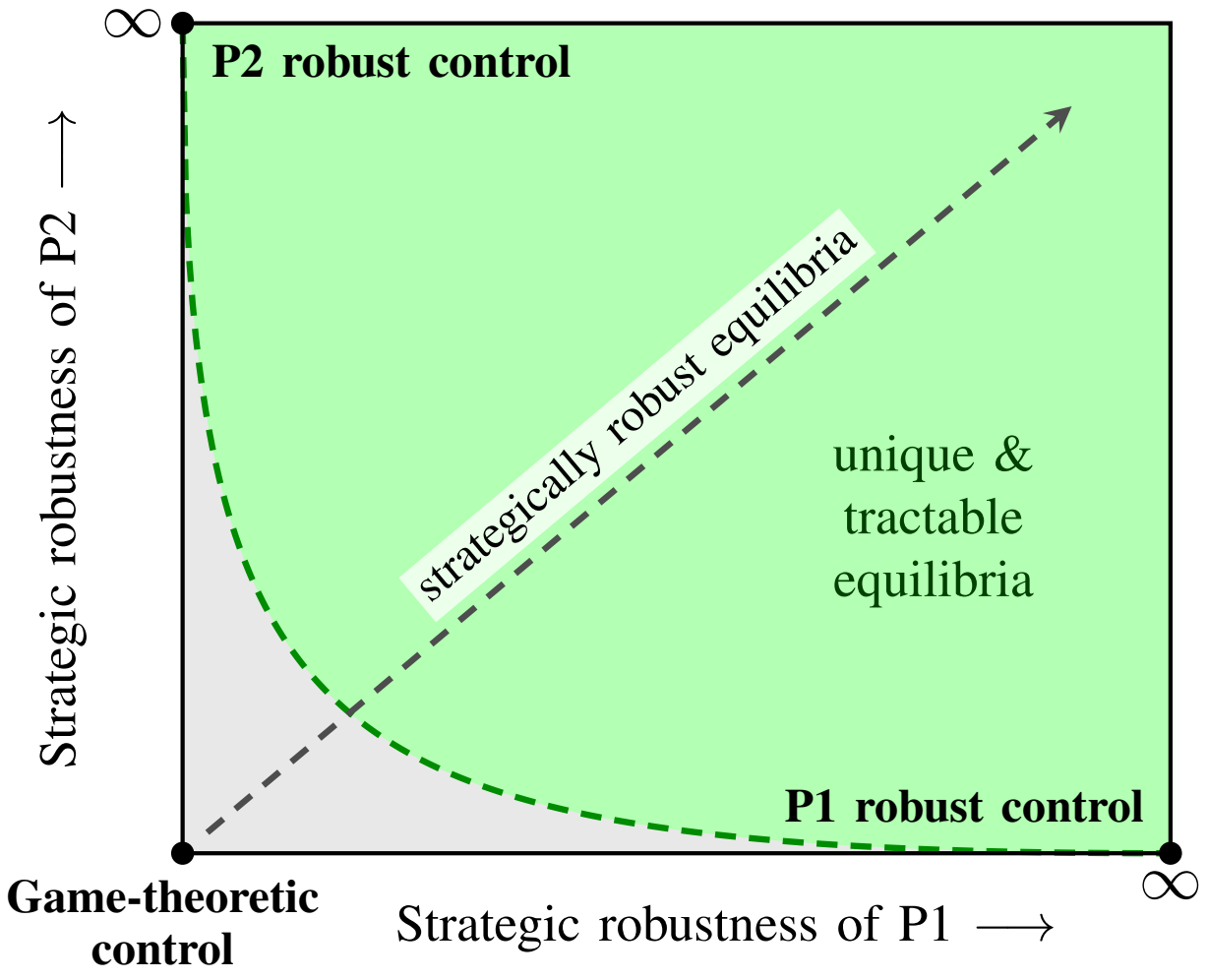

In game-theoretic control, agents are modeled as rational decision-makers minimizing their own control cost.
The natural solution concept is that of Nash equilibrium, a configuration at which no agent can lower their cost by unilaterally changing their control policy.
While natural in theory, the practical deployment of Nash equilibrium control policies suffers from a fundamental limitation: policies are brittle, i.e., sensitive to misspecification of the other agents’ control policies. For instance, \cite{nagel1995unraveling} demonstrates that adhering to a standard Nash equilibrium when opponents are non-rational can lead to significantly suboptimal outcomes.

In addition to their poor robustness, Nash equilibria also pose significant computational challenges, and can generally be computed efficiently only under structural assumptions, such as monotonicity of the pseudo-gradient or the invertibility conditions arising in linear-quadratic games~\cite{bacsar1998dynamic}, which restrict the class of problems that can be solved in practice.

In robust control, instead, the other agents are modeled as adversarial disturbances, whose goal is to maximize the control cost.
The classical $\mathcal{H}_\infty$ control problem can be cast as a zero-sum game between the control engineer, who minimizes control cost, and a fictitious agent, who instead maximizes control cost and whose influence is attenuated only through a penalty on its magnitude.
A natural way to inject robustness to strategic uncertainty is thus to treat the control policies of the others as exogenous disturbances, as in, e.g.,~\cite{amato2002guaranteeing,jimenez2006varepsilon,pantazis2023data}.
Doing so, however, completely disregards the underlying game-theoretic interactions: the other agents are not arbitrary disturbances, but decision-makers with their own objectives, whose behavior should be expected to remain close to the policies prescribed by the game.
Thus, while computationally tractable, it inevitably leads to overly conservative behavior and sacrifices nominal performance. In other words, game-theoretic control trusts the other agents entirely, while robust control treats their control inputs as arbitrary disturbances.

In this paper, we populate the spectrum between these two solution concepts; see Figure~\ref{fig:meta}. We take inspiration from recent work on strategic risk aversion in multi-agent reinforcement learning~\cite{mazumdar2025tractable} and strategically robust game theory in one-shot games~\cite{lanzetti2025strategically}, from which we also inherit the name ``strategic robustness.''
Specifically, we propose that agents choose their strategy not by playing against the policies of the other agents (as in game-theoretic control) but instead against a fictitious adversary who seeks to inflict maximum damage (as in robust control), but pays a penalty for deviating from the policies of the other agents.
This deviation penalty controls how tightly the adversary is anchored to the other players' policies.
When it is prohibitively large, deviations are too expensive and the adversary is forced toward the policies of the other agents, recovering standard game-theoretic control.
As the penalty weakens, the adversary gains power and the player's policy becomes increasingly robust. When an explicit disturbance penalty keeps the adversarial problem well posed, the zero-deviation-penalty limit is a classical worst-case $\mathcal{H}_\infty$ control problem---a connection we make precise in this paper.
Together with the direct disturbance penalty, the deviation penalty defines a two-parameter family containing the standard Nash and $\mathcal{H}_\infty$ problems as distinct limiting cases.

Beyond our formulation, our work also highlights other advantages of strategically robust (game-theoretic) control. In particular, it often enjoys better computational properties compared to standard game-theoretic control and can sometimes lead to more efficient and performant outcomes.

We instantiate these ideas in the setting of discrete-time dynamic games. Dynamic games are used in a wide variety of control applications, including robotics~\cite{fridovich2020efficient}, aerospace~\cite{chai2020linear}, and economics~\cite{pindyck2003optimal}, to name a few.
In these games, the state of a dynamical system is controlled by a set of agents that select policies to minimize their personal control cost, which depends on the state of the system and on the agents' inputs.
Specifically, we study two classes of dynamic games: the finite-horizon open-loop setting, with linear dynamics, convex costs, and convex constraints, and the finite-horizon closed-loop linear-quadratic setting, which we further extend to the infinite horizon.

\subsubsection*{Contributions}
We introduce strategically robust dynamic equilibria, a solution concept for dynamic games that protects agents against strategic uncertainty, i.e., uncertainty about the behavior of other agents.
Strategically robust dynamic equilibria allow us to interpolate between game-theoretic control (Nash equilibria) and robust control ($\mathcal H_\infty$), all of this while being ``more computationally tractable'' than standard Nash equilibria.
Specifically, our contributions are fourfold.
\begin{enumerate}
    \item \emph{Open loop.} After establishing existence of strategically robust dynamic equilibria, we give verifiable conditions under which a reweighted pseudo-gradient is strongly monotone.
    These conditions imply a unique equilibrium and global linear convergence of the projected pseudo-gradient method on compact feasible sets, or when the pseudo-gradient is globally Lipschitz.

    \item \emph{Closed loop, finite horizon.}
    We characterize closed-loop strategically robust dynamic equilibria via coupled backward Riccati equations. This way, we show that, under standard conditions, the strategically robust dynamic equilibrium is unique (in the subgame-perfect sense), linear, and Markovian. Strategic robustness additionally yields a sufficient condition for invertibility of the coupled stage system.

    \item \emph{Closed loop, infinite horizon.} Outside of the zero-sum setting, Nash equilibria of infinite-horizon games are notoriously hard to compute.
    We show that, with sufficient robustness, a stationary linear equilibrium exists, is unique in a ball that becomes arbitrarily large as agents become more robust, and is reached at a geometric rate by best-response iteration---every step of which is a classical $\mathcal{H}_\infty$ design.

    \item \emph{Experiments.} Strategic robustness protects agents against misspecified policies and adversarial perturbations of other agents' controllers. In our examples, we empirically demonstrate that it can also lead to a ``free-lunch'' effect, in that it coordinates the agents and improves individual nominal performance or social welfare.
\end{enumerate}

\subsubsection*{Related work}
Unlike the strategically robust equilibria of~\cite{lanzetti2025strategically}, introduced for static games (and instantiated in the aggregative setting in~\cite{feik2026strategically} and in multi-agent trajectory optimization in~\cite{qin2026strategically}), we operate in a dynamic setting and seek control policies rather than static decisions, using deterministic rather than mixed policies and replacing the hard-constrained ambiguity set of~\cite{lanzetti2025strategically} by a deviation penalty.
Compared with risk-averse quantal equilibria~\cite{mazumdar2025tractable,zhang2026provably,zhang2025convergent}, we do not model bounded rationality and instead consider a standard control-theoretic setting with continuous state and action spaces. The augmented-game reformulation that we use to derive our results is, nonetheless, inspired by that line of work.
Similarly to~\cite{lanzetti2025strategically,feik2026strategically,qu2026training}, we observe empirically that strategic risk aversion and strategic robustness can both protect and coordinate agents.
The closely related work~\cite{rabbani2026optimal} studies input disturbances in two-player linear-quadratic games. There, player 1 compensates online for the observed disturbance, while player 2 does not best-respond; the resulting law is, as the authors note, not itself an equilibrium. Robustness is therefore reactive rather than strategic, ``on top of'' the standard Nash equilibrium.
More broadly, related work on robust dynamic games~\cite{jank2002disturbed,vandenbroek2003robust,jimenez2020robust,roy2025soft} studies games where the dynamics of each player are affected by an exogenous disturbance, as in standard robust control. In contrast, our perspective treats uncertainty endogenously. This leads to a fundamentally different formulation which we study in this paper and use to interpolate between game-theoretic control and robust control. Nevertheless, in the linear-quadratic setting with specific parameter choices, our formulation can recover that with exogenous uncertainty.
Finally, a preliminary version of this work, confined to the finite-horizon closed-loop linear-quadratic setting, appeared in~\cite{velasevic2026strategically}.
The present paper subsumes and substantially extends it. First, we
treat the open-loop setting with general convex costs and
constraints, where we establish existence and give conditions for strong monotonicity (Section~\ref{sec:open_loop}). Second, we make the interpolation
between Nash equilibria and $\mathcal{H}_\infty$ control precise,
exhibiting both as limits of a two-parameter family (Section~\ref{sec:hinf}).
Third, we study the infinite horizon, proving existence, uniqueness, and geometric convergence of best-response iteration in a neighborhood of the decoupled
$\mathcal{H}_\infty$ gains (Section~\ref{sec:infinite}). Finally, we include open-loop experiments on platooning and capacity accumulation in Section~\ref{sec:experiments}.

\subsubsection*{Notation}
For a symmetric matrix $Q$, we write $Q\succeq 0$ ($Q\succ 0$) if $Q$ is positive semidefinite (definite), and $\norm{x}_Q^2 \coloneqq x^\top Q x$. The notation $\norm{\cdot}$ denotes the Euclidean norm for vectors and the induced spectral norm for matrices; for derivatives between matrix spaces, it denotes the corresponding operator norm. We denote by $\lambda_{\min}(M)$ and $\lambda_{\max}(M)$ the smallest and largest eigenvalues of a symmetric matrix $M$, and by $\sigma_{\min}(M) = \sqrt{\lambda_{\min}(M^\top M)}$ and $\sigma_{\max}(M) = \sqrt{\lambda_{\max}(M^\top M)}$ the smallest and largest singular values. For player $i$, the superscript $-i$ collects the corresponding quantities of all other players. The matrix $I$ denotes the identity of the appropriate dimension.

\section{Dynamics, costs, and strategic robustness}\label{sec:setting}

We consider a discrete-time $N$-player dynamic game over a finite horizon $T\in\naturals_{>0}$.
The system state $x_t \in \reals^n$ evolves according to
\begin{align}
    x_{t+1} &= A_t x_t + \sum_{j=1}^N B_t^j u_t^j \label{eq:system_dynamics}
\end{align}
where $u_t^i \in \reals^{n_i}$ is the control input of player $i$ which satisfies $(u_0^i, u_1^i, \ldots, u_{T-1}^i) \in \mathcal{U}^i$, $A_t$ and $B_t^i$ are matrices of appropriate dimensions, $x_0\in\reals^n$ is a given initial condition, and the constraint set $\mathcal{U}^i$ is the set that the control sequence is restricted to.
We consider two cases of the information structure:
\begin{enumerate}
    \item \emph{Open loop}: The input sequence $u_t^i$ is committed to at time $0$ and does not react to the realized state.
    \item \emph{Closed loop}: The input sequence is generated by causal policies, which are measurable functions of the history.
\end{enumerate}
In the remainder of the section, we denote both by $\pi_t^i$ for brevity---in the open-loop case we set $u_t^i \equiv \pi_t^i$, in the closed-loop case $u_t^i \equiv \pi_t^i(x_0, \ldots, x_t)$.
The corresponding action sets are defined in Sections~\ref{sec:open_loop} and~\ref{sec:closed_loop}. In either case, each player $i$ seeks to minimize the cost functional
\begin{equation}\label{eq:cost}
    l_T^i(x_T)
    +
    \sum_{t=0}^{T-1} l_t^i(x_t, u_t^i, u_t^{-i}),
\end{equation}
subject to the system dynamics and constraints~\eqref{eq:system_dynamics}.

Write $J_{\mathrm{nom}}^i(\pi^i,\pi^{-i})$ for~\eqref{eq:cost} evaluated along
the trajectory generated by $(\pi^i,\pi^{-i})$. A profile
$(\bar\pi^1,\ldots,\bar\pi^N)$ is a \emph{Nash equilibrium} if
$J_{\mathrm{nom}}^i(\bar\pi^i,\bar\pi^{-i})\le
J_{\mathrm{nom}}^i(\pi^i,\bar\pi^{-i})$ for every player $i$ and every
admissible $\pi^i$: no player lowers their cost by deviating alone.
Unfortunately, Nash equilibria are brittle, and their performance can quickly deteriorate when the other players deviate from Nash equilibrium feedback policies.

\begin{example}[motivating]\label{example:motivating}
    Consider a two-player closed-loop setting with scalar dynamics $x_{t+1} = 1.05 x_t + u^1_t + u^2_t$ and control costs $\sum_{t = 0}^{T - 1} (u^i_t)^2 + Q_T^i x_T^2$.
    Under the Nash equilibrium closed-loop policies, the state evolves as $x_{t + 1} \approx 0.98x_t$ for the majority of the horizon. As the coefficient is close to 1, small deviations in agents' policies can drastically impact the cost.
\end{example}

This fragility prompts players to seek protection, so we consider the strategically robust control cost
\begin{align}
    J^i(\pi^i,\pi^{-i}) = \sup_{\sigma^i} & \quad l_T^i(x_T) \nonumber
    \\ &+
    \sum_{t=0}^{T-1} l_t^i(x_t, u_t^i, d_t^i) - p_t^i(d_t^i, u_t^{-i}). \label{eq:robust_cost}
\end{align}
Let $x_t$ denote the fictitious trajectory induced by $(\pi^i,\sigma^i)$: $x_{t+1}=A_tx_t+B_t^iu_t^i+B_t^{-i}d_t^i$. The supremum is taken over the policies $\sigma^i$ of the fictitious adversary. The action vector $d_t^i\in\reals^{\sum_{j\neq i}n_j}$ is generated by $\sigma^i$, and $B_t^{-i}$ is formed by horizontally concatenating the corresponding matrices $B_t^j$. The function $p_t^i : \reals^{\sum_{j \neq i} n_j} \times \reals^{\sum_{j \neq i} n_j} \to \reals$ penalizes the deviation of the adversary from other players. We assume $p_t^i(y,y)=0$, $p_t^i(x,y)\ge0$, and take $p_t^i$ to be convex and differentiable in both arguments.  Different players generally evaluate their robust costs along different fictitious trajectories.
In other words, player $i$ plays against a fictitious adversary that maximizes their control cost while paying $p_t^i(d_t^i,u_t^{-i})$ for deviating from the policies $\pi^{-i}$.

This cost leads to a notion of equilibrium in which players are robust to each other's policies---the dynamic extension of the strategically robust equilibria of~\cite{lanzetti2025strategically}---made precise in Definitions~\ref{def:sre:ol} and~\ref{def:sre:cl}.

\begin{remark}[level of robustness]\label{rem:level_of_robustness}
We often take $p_t^i(d_t^i, u_t^{-i}) =  \norm{d_t^i - u_t^{-i}}^2_{M_t^i}$, where a common choice is $M_t^i=m_t^iI$ with $m_t^i>0$, specifying the level of robustness. Specifically:
\begin{itemize}
    \item No robustness: as $M_t^i \to \infty$, deviations become prohibitively expensive, the adversary is forced toward $u_t^{-i}$, and each player faces their nominal best-response problem; the limiting equilibrium is the Nash equilibrium of the nominal game.

    \item Robust regime: for finite $M_t^i$ the adversary can explore deviations from $u_t^{-i}$, preparing player $i$ for misspecified or adversarial behavior, and smaller $M_t^i$ increases robustness. The limit $M_t^i\to0$ removes the deviation penalty; when a direct penalty on $d_t^i$ keeps the adversarial problem well posed, this is the $\mathcal{H}_\infty$ limit of Section~\ref{sec:hinf}.

    \item Heterogeneous robustness: choosing $M_t^i$ block diagonal, $M_t^i = \text{diag}(M_t^{i,j})_{j \neq i}$, lets player $i$ parameterize their robustness to each opponent $j \neq i$ separately, modeling players perceived as more reliable than others.
\end{itemize}
\end{remark}

Let $l^i(x,u^i,d^i)\coloneqq l_T^i(x_T)+\sum_{t=0}^{T-1}l_t^i(x_t,u_t^i,d_t^i)$ and $p^i(d^i,u^{-i})\coloneqq\sum_{t=0}^{T-1}p_t^i(d_t^i,u_t^{-i})$. These are
functions of a state trajectory and of input sequences, and the arguments only name the slots; what is substituted into them is prescribed by the policies.
Let $x(\pi^i,\sigma^i)$ be the fictitious trajectory generated from $x_0$ by $\pi^i$ and $\sigma^i$, and let $u^i$, $d^i$, and $u^{-i}$ be the input sequences that $\pi^i$, $\sigma^i$, and $\pi^{-i}$ produce along it. The \emph{augmented cost} of player $i$ is
$J_{\mathrm{aug}}^i(\pi^i,\sigma^i,\pi^{-i})\coloneqq
l^i(x(\pi^i,\sigma^i),u^i,d^i)-p^i(d^i,u^{-i})$, and the robust
cost~\eqref{eq:robust_cost} is $J^i=\sup_{\sigma^i}J_{\mathrm{aug}}^i$. In the
open-loop setting policies are themselves input sequences, so we write
$J_{\mathrm{aug}}^i(u^i,d^i,u^{-i})$ there.
We work under one of two sets of assumptions, according to the setting.

\begin{assumption}[convex open-loop]\label{ass:ol}
For every player $i$, the induced feasible sequence set $\mathcal U^i$ is nonempty, closed, and convex; $l_t^i$
and $p_t^i$ are twice continuously differentiable; and, at every joint profile,
the blocks of the Hessians of $l_t^i$ for all $t \in \{0, \ldots, T-1\}$ obey
\begin{equation}\label{eq:ass:curvature}
\begin{gathered}
    \underline q^iI\preceq\nabla^2_{xx}l_t^i\preceq\overline q^iI,
    \quad
    \nabla^2_{u^iu^i}l_t^i\succeq\underline r^iI,
    \quad
    \nabla^2_{d^id^i}l_t^i\preceq\overline\ell^iI,\\
    \norm{\nabla^2_{xu^i}l_t^i}\le\kappa_u^i,
    \qquad
    \norm{\nabla^2_{xd^i}l_t^i}\le\kappa_d^i,
\end{gathered}
\end{equation}
and $\underline q^iI\preceq\nabla^2_{xx}l_T^i\preceq\overline q^iI$, with $\underline q^i\ge0$ and $\underline r^i>0$, while the penalty obeys
$\nabla^2_{d^id^i}p_t^i\succeq\mu_p^iI$ with $\mu_p^i>0$ and
$\norm{\nabla^2_{d^iu^j}p_t^i}\le c^{ij}$ for every $j\neq i$.
\end{assumption}

The own-input block is bounded below and the adversary block above, with
$\overline\ell^i$ free to be negative when the cost penalizes the adversary
directly. The cross blocks are bounded in norm only and vanish for stage costs
separable in $(x,u^i,d^i)$; Section~\ref{sec:open_loop} combines all of these
into the curvature constants~\eqref{eq:constants}. The constants $c^{ij}$
measure the only channel through which distinct players interact.

\begin{assumption}[LQ closed-loop]\label{ass:cl}
For every player $i$, $\mathcal{U}^i= \reals^{n_i} \times \cdots \times \reals^{n_i}$,
$l_t^i(x_t,u_t^i,u_t^{-i})=\norm{x_t}^2_{Q_t^i}+\norm{u_t^i}^2_{R_t^i}-\norm{u_t^{-i}}^2_{S_t^i}$,
and $p_t^i(d_t^i,u_t^{-i})=\norm{d_t^i-u_t^{-i}}^2_{M_t^i}$, with $Q_t^i,Q_T^i\succeq0$, $R_t^i\succ0$, $S_t^i\succeq0$, and $M_t^i\succ0$.
\end{assumption}

Assumption~\ref{ass:ol} is in force throughout Section~\ref{sec:open_loop}, Assumption~\ref{ass:cl} throughout Section~\ref{sec:closed_loop}.

\section{The open-loop setting: existence, monotonicity, and computation}\label{sec:open_loop}
Throughout this section Assumption~\ref{ass:ol} is in force. Players choose
input sequences directly, so a policy is the sequence itself, $\pi^i\equiv u^i$.
The action sets are $
\mathcal U^i,
\mathcal U^{-i}\coloneqq\prod_{j\neq i}\mathcal U^j,
$
where the first is the set of all sequences that satisfy
constraints. An adversarial policy is a deviation
sequence $\sigma^i\equiv d^i\in\mathcal U^{-i}$. These sets are nonempty, closed,
and convex.

To simplify notation, we let
$x=\operatorname{col}(x_0,\ldots,x_T)$,
$u^i=\operatorname{col}(u_0^i,\ldots,u_{T-1}^i)$, and
$d^i=\operatorname{col}(d_0^i,\ldots,d_{T-1}^i)$, and write
$\Phi(t,s)\coloneqq A_{t-1}\cdots A_s$ for $t>s$, with $\Phi(t,t)\coloneqq I$.
Iterating the fictitious dynamics gives $x=G^iu^i+G^{-i}d^i+Cx_0$, where
\begin{equation*}
    {\setlength{\arraycolsep}{2.5pt}
    G^i=
    \begin{bmatrix}
        0 & 0 & \cdots & 0\\
        B_0^i & 0 & \cdots & 0\\
        A_1B_0^i & B_1^i & \cdots & 0\\
        \vdots & \vdots & \ddots & \vdots\\
        \Phi(T,1)B_0^i & \Phi(T,2)B_1^i & \cdots & B_{T-1}^i
    \end{bmatrix},
    C = \begin{bmatrix} I \\ A_0 \\ \Phi(2,0) \\ \vdots \\ \Phi(T,0) \end{bmatrix},}
\end{equation*}
and $G^{-i}$ is
$G^i$ with every $B_s^i$ replaced by $B_s^{-i}$. Substituting this into $l^i$
turns the dynamic game into a static, finite-dimensional problem in $(u,d)$.

Elimination of the state turns the stagewise bounds of Assumption~\ref{ass:ol} into two
curvature constants per player,
\begin{align}\label{eq:constants}
    \mu_U^i&\coloneqq\underline r^i+\underline q^i\sigma_{\min}(G^i)^2
        -2\kappa_u^i\sigma_{\max}(G^i),\nonumber\\
    \mu_D^i&\coloneqq\mu_p^i-\overline\ell^i-\overline q^i\sigma_{\max}(G^{-i})^2
        -2\kappa_d^i\sigma_{\max}(G^{-i}),
\end{align}
which the chain-rule bound from Lemma~\ref{lemma:transfer} shows to satisfy
$\nabla^2_{u^iu^i}J^i_{\mathrm{aug}}\succeq\mu_U^iI$ and
$-\nabla^2_{d^id^i}J^i_{\mathrm{aug}}\succeq\mu_D^iI$. So $\mu_U^i$
lower-bounds the curvature of player $i$'s problem in their own input and
$\mu_D^i$ that of their adversary's, while the constants $c^{ij}$ of
Assumption~\ref{ass:ol} bound the interference between distinct players.
Every result in this section asks for $\mu_U^i>0$ and $\mu_D^i>0$. It makes each $J^i_{\mathrm{aug}}$ strongly convex in $u^i$ and strongly concave in $d^i$, hence coercive in both
(Section~\ref{sec:proofs_ol}). This requirement is satisfied in many common cost structures. Losses that are decoupled in $(x, u)$ make $\kappa_u^i = \kappa_d^i = 0$, so $\mu_U^i > 0$ is always satisfied and $\mu_D^i > 0$ can be satisfied by choosing a sufficiently convex penalty.

\subsection{Equilibria and the augmented game}\label{sec:ol_static}

We now explicitly define the open-loop strategically robust dynamic equilibrium.

\begin{definition}[open-loop strategically robust dynamic equilibrium]\label{def:sre:ol}
    A tuple $(\bar u^{1},\ldots,\bar u^N)\in\mathcal U^1\times\cdots\times\mathcal U^N$ forms an open-loop strategically robust dynamic equilibrium if for all players
    \begin{equation}\label{eq:sre:ol}
        J^i(\bar u^i,\bar u^{-i})
        \leq
        J^i(u^i,\bar u^{-i})
        \qquad
        \forall u^i\in\mathcal U^i.
    \end{equation}
\end{definition}

Definition~\ref{def:sre:ol} reads as a two-stage game: the players commit to $u^i$, and the fictitious adversaries then respond. Under the conditions below it is equivalent to an augmented game in which all of them move at once.

\begin{definition}[equilibrium of the open-loop augmented game]\label{def:aug:ol}
    The sequences $(\bar u^1,\bar d^{1}, \ldots, \bar u^N,\bar d^{N})$ are an equilibrium of the augmented game if for all players $i$ we have
    \begin{equation}\label{eq:aug:ol}
        J_\mathrm{aug}^i(\bar u^{i}, d^i, \bar u^{-i})
        \leq
        J_\mathrm{aug}^i(\bar u^{i}, \bar d^i, \bar u^{-i})
        \leq
        J_\mathrm{aug}^i(u^i, \bar d^i, \bar u^{-i})
    \end{equation}
    for all $u^i\in\mathcal U^i$ and $d^i\in\mathcal U^{-i}$.
\end{definition}

Under Assumption~\ref{ass:ol} with $\mu_U^i,\mu_D^i>0$, the two formulations are equivalent.

\begin{lemma}[equivalence in the open-loop setting]\label{lemma:equivalence:ol}\,
\begin{enumerate}
    \item [\textnormal{(a)}]
    \emph {Necessity.} If
    $(\bar u^1,\bar d^{1},\ldots,\bar u^N,\bar d^{N})$ is an equilibrium of the
    augmented game, then $(\bar u^1,\ldots,\bar u^N)$ is an open-loop
    strategically robust dynamic equilibrium and each $\bar d^i$ is a worst-case
    sequence in~\eqref{eq:robust_cost}, i.e., a best response to $\bar u^i$.
    \item [\textnormal{(b)}]
    \emph{Sufficiency.} Assume that $\mu_U^i > 0, \mu_D^i > 0$ for all $i$. Every open-loop strategically robust dynamic equilibrium
    $(\bar u^1,\ldots,\bar u^N)$ admits worst-case sequences
    $\bar d^1,\ldots,\bar d^N$ such that
    $(\bar u^1,\bar d^{1},\ldots,\bar u^N,\bar d^{N})$ is an equilibrium of the
    augmented game.
\end{enumerate}
\end{lemma}

Lemma~\ref{lemma:equivalence:ol} is proved in Section~\ref{sec:proofs_ol}, by specializing the general equivalence principle of Lemma~\ref{lemma:equivalence:gen} and verifying its saddle-point hypothesis under the assumptions above.

\subsection{Existence}
The following is a standard proposition that guarantees existence in the open-loop case.

\begin{proposition}[existence in the open-loop setting]\label{prop:ole}
    Let Assumption~\ref{ass:ol} hold with $\mu_U^i,\mu_D^i>0$ for every $i$.
    \begin{enumerate}
        \item[\textnormal{(a)}] If every $\mathcal U^i$ is compact, a strategically
        robust dynamic equilibrium $({u^1}^\star,\ldots,{u^N}^\star)$ exists, together
        with worst-case adversarial sequences $({d^1}^\star,\ldots,{d^N}^\star)$.
        \item[\textnormal{(b)}] If $\mathcal U^i=\reals^{Tn_i}$ for every $i$,
        then every solution of the stationarity system
        \begin{equation}\label{eq:ol:foc}
            \nabla_{u^i}J^i_{\mathrm{aug}}=0,
            \qquad
            \nabla_{d^i}J^i_{\mathrm{aug}}=0,
            \qquad i=1,\ldots,N,
        \end{equation}
        is a strategically robust dynamic equilibrium.
    \end{enumerate}
\end{proposition}

The conditions $\mu_U^i,\mu_D^i>0$ are sufficient and not necessary, as the proof uses only strict convexity--concavity and coercivity. Part (b) allows us to attempt to find such equilibria in the unconstrained case using various solvers; however, these might fail as they are sensitive to initial conditions.

\subsection{Monotonicity analysis}\label{sec:monotonicity}
Beyond existence, strong monotonicity ensures uniqueness of the equilibrium and, on compact sets, global linear convergence of the projected pseudo-gradient method for a sufficiently small fixed step size; see, e.g.,~\cite{facchinei2003finite}.

For $z\coloneqq\operatorname{col}(u^1,\ldots,u^N,d^1,\ldots,d^N)$, define the pseudo-gradient
$
F(z)\coloneqq\operatorname{col}\big((\nabla_{u^i}J_{\mathrm{aug}}^i)_{i=1}^N,\,(-\nabla_{d^i}J_{\mathrm{aug}}^i)_{i=1}^N\big).
$
For $\mu>0$, we say that $F$ is $\mu$-strongly monotone if
$\langle F(z)-F(z'),z-z'\rangle\geq\mu\norm{z-z'}^2$ for all $z,z'$. Since $F$ is $C^1$ under Assumption~\ref{ass:ol}, this is equivalent to
$\tfrac12(\nabla F(z)+\nabla F(z)^\top)\succeq\mu I$ for all $z$; we use this characterization throughout.

More generally, for weights $\lambda=(\lambda_1,\ldots,\lambda_N)\in\reals_{>0}^N$, let
$
F_\lambda(z)\coloneqq\operatorname{col}\big((\lambda_i\nabla_{u^i}J_{\mathrm{aug}}^i)_{i=1}^N,\,(-\lambda_i\nabla_{d^i}J_{\mathrm{aug}}^i)_{i=1}^N\big).
$
We call the game strongly $\lambda$-monotone if $F_\lambda$ is strongly monotone; equal weights recover the standard notion. Rescaling each player's loss leaves the saddle inequalities defining the augmented-game equilibrium unchanged, so such a certificate gives the same equilibrium while balancing the players' curvatures and couplings.

For each player $k$, let
\begin{equation}\label{eq:cycle_constants}
    c_{\mathrm{in}}^k\coloneqq\sum_{i\neq k}c^{ik},
    \quad
    c_{\mathrm{out}}^k\coloneqq\sum_{j\neq k}c^{kj},
    \quad
    \chi^k\coloneqq\frac{4\mu_U^k\mu_D^k}{c_{\mathrm{in}}^k\,c_{\mathrm{out}}^k},
\end{equation}
Let $\chi^{(1)}\le \chi^{(2)}$ be the two smallest of $\chi^1,\ldots,\chi^N$, and let player~$k^\star$ attain $\chi^{(1)}$.

\begin{theorem}[certificate of strong $\lambda$-monotonicity]\label{thm:monotonicity}
    Let Assumption~\ref{ass:ol} hold with $\mu_U^k,\mu_D^k>0$ for every $k$. If
    \begin{equation}\label{eq:monotonicity_condition}
        \chi^{(1)}\chi^{(2)}>1,
    \end{equation}
    then the game is strongly $\lambda$-monotone for the weights
    \begin{equation}\label{eq:weights}
        \lambda_k=\frac{c_{\mathrm{in}}^k}{\mu_U^k}\ \ (k\neq k^\star),
        \qquad
        \lambda_{k^\star}=\sqrt{\frac{\chi^{(1)}}{\chi^{(2)}}}\,\frac{c_{\mathrm{in}}^{k^\star}}{\mu_U^{k^\star}}.
    \end{equation}
\end{theorem}

Both sums in~\eqref{eq:cycle_constants} are positive under Assumption~\ref{ass:ol}, and the ratio $\chi^k$ compares the curvature of player~$k$ and their adversary with the total coupling entering and leaving them; only the two smallest ratios enter the condition. To demonstrate that these conditions can be relatively mild, for two players and in the LQ case, the theorem specializes as follows.

\begin{corollary}[two-player LQ monotonicity]\label{cor:two_player_monotonicity}
Consider the open-loop LQ game with $N=2$ and costs as in Assumption~\ref{ass:cl}, with $Q_t^i=Q_T^i=Q^i$, $R_t^i=R^i$, $S_t^i=s^iI$, and $M_t^i=m^iI$. Set
$
    a^i\coloneqq\lambda_{\min}(R^i)+\lambda_{\min}(Q^i)\sigma_{\min}(G^i)^2,
    b^i\coloneqq m^i+s^i-\lambda_{\max}(Q^i)\sigma_{\max}(G^{-i})^2.
$
If
\begin{equation}\label{eq:two_player_lq}
    b^1,b^2>0,
    \qquad
    (m^1m^2)^2<16\,a^1b^1a^2b^2,
\end{equation}
the game is strongly $\lambda$-monotone for the weights~\eqref{eq:weights}.
\end{corollary}

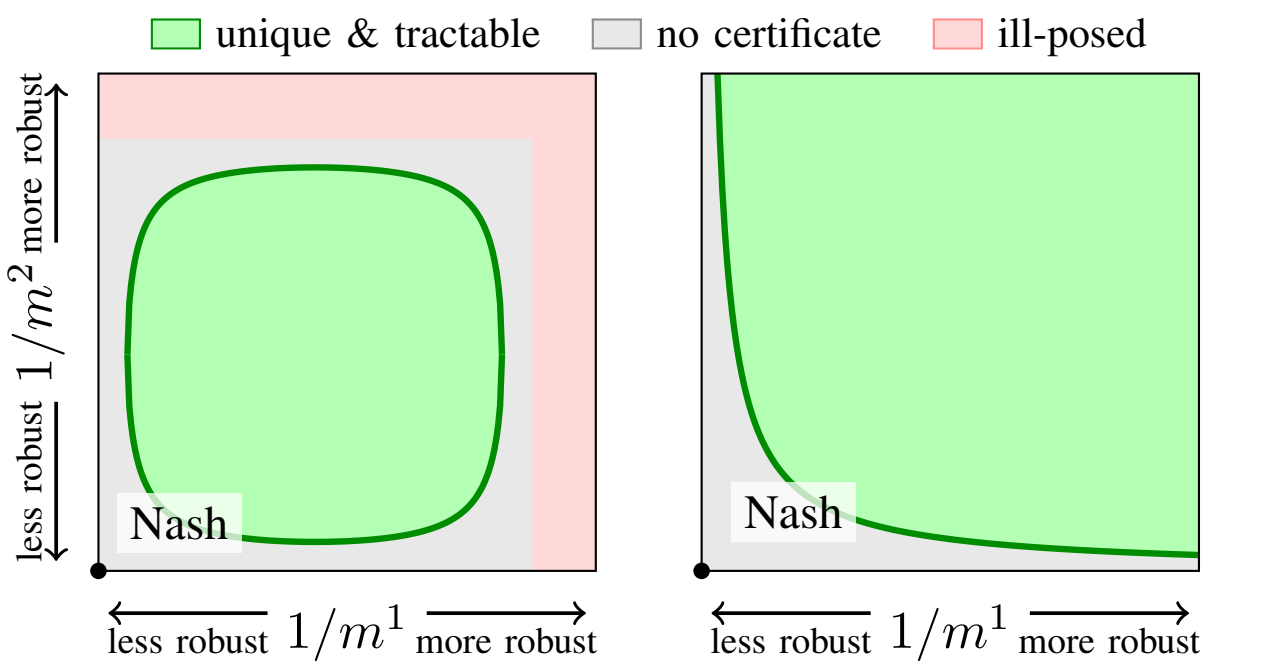
\begin{figure}[!t]
    \centering
    \begin{tikzpicture}
    \begin{groupplot}[
        group style={group size=2 by 1, horizontal sep=0.7cm, ylabels at=edge left},
        width=0.55\columnwidth, height=0.55\columnwidth,
        xlabel={$1/m^1$}, ylabel={$1/m^2$},
        xtick=\empty, ytick=\empty,
        axis on top, enlargelimits=false, clip=false,
        xlabel style={at={(axis description cs:0.5,-0.03)}, anchor=north},
        ylabel style={at={(axis description cs:-0.03,0.5)}, anchor=south},
        /tikz/robarrow/.style={->, line width=0.7pt},
        /tikz/roblabel/.style={font=\scriptsize, inner sep=2.5pt},
    ]
    \nextgroupplot[xmin=0, xmax=1.15, ymin=0, ymax=1.15]
        \addplot[fill=red!15, draw=none]
            coordinates {(0,0) (1.15,0) (1.15,1.15) (0,1.15)} \closedcycle;
        \addplot[fill=gray!18, draw=none]
            coordinates {(0,0) (1,0) (1,1) (0,1)} \closedcycle;
        \addplot[fill=green!30, draw=none] coordinates {(0.0670,0.5066) (0.0862,0.7274) (0.1055,0.7905) (0.1247,0.8269) (0.1440,0.8510) (0.1632,0.8682) (0.1825,0.8811) (0.2017,0.8911) (0.2210,0.8990) (0.2402,0.9054) (0.2594,0.9107) (0.2787,0.9151) (0.2979,0.9187) (0.3172,0.9217) (0.3364,0.9243) (0.3557,0.9264) (0.3749,0.9282) (0.3942,0.9296) (0.4134,0.9308) (0.4326,0.9317) (0.4519,0.9323) (0.4711,0.9328) (0.4904,0.9330) (0.5096,0.9330) (0.5289,0.9328) (0.5481,0.9323) (0.5674,0.9317) (0.5866,0.9308) (0.6058,0.9296) (0.6251,0.9282) (0.6443,0.9264) (0.6636,0.9243) (0.6828,0.9217) (0.7021,0.9187) (0.7213,0.9151) (0.7406,0.9107) (0.7598,0.9054) (0.7790,0.8990) (0.7983,0.8911) (0.8175,0.8811) (0.8368,0.8682) (0.8560,0.8510) (0.8753,0.8269) (0.8945,0.7905) (0.9138,0.7274) (0.9330,0.5066) (0.9330,0.4934) (0.9138,0.2726) (0.8945,0.2095) (0.8753,0.1731) (0.8560,0.1490) (0.8368,0.1318) (0.8175,0.1189) (0.7983,0.1089) (0.7790,0.1010) (0.7598,0.0946) (0.7406,0.0893) (0.7213,0.0849) (0.7021,0.0813) (0.6828,0.0783) (0.6636,0.0757) (0.6443,0.0736) (0.6251,0.0718) (0.6058,0.0704) (0.5866,0.0692) (0.5674,0.0683) (0.5481,0.0677) (0.5289,0.0672) (0.5096,0.0670) (0.4904,0.0670) (0.4711,0.0672) (0.4519,0.0677) (0.4326,0.0683) (0.4134,0.0692) (0.3942,0.0704) (0.3749,0.0718) (0.3557,0.0736) (0.3364,0.0757) (0.3172,0.0783) (0.2979,0.0813) (0.2787,0.0849) (0.2594,0.0893) (0.2402,0.0946) (0.2210,0.1010) (0.2017,0.1089) (0.1825,0.1189) (0.1632,0.1318) (0.1440,0.1490) (0.1247,0.1731) (0.1055,0.2095) (0.0862,0.2726) (0.0670,0.4934)} \closedcycle;
        \addplot[green!55!black, very thick, domain=0.0669:0.9331, samples=200]
            {0.5 + 0.5*sqrt(max(0, 1 - 1/(16*x*(1-x))))};
        \addplot[green!55!black, very thick, domain=0.0669:0.9331, samples=200]
            {0.5 - 0.5*sqrt(max(0, 1 - 1/(16*x*(1-x))))};
        \node[circle, fill=black, inner sep=1.1pt] at (axis cs:0,0) {};
        \node[anchor=south west, font=\small, inner sep=2.5pt,
              fill=white, fill opacity=0.75, text opacity=1]
              at (axis cs:0.04,0.04) {Nash};
        \draw[robarrow] (axis description cs:0.34,-0.085) -- (axis description cs:0.02,-0.085);
        \node[roblabel, anchor=north] at (axis description cs:0.18,-0.085) {less robust};
        \draw[robarrow] (axis description cs:0.66,-0.085) -- (axis description cs:0.98,-0.085);
        \node[roblabel, anchor=north] at (axis description cs:0.82,-0.085) {more robust};
        \draw[robarrow] (axis description cs:-0.085,0.34) -- (axis description cs:-0.085,0.02);
        \node[roblabel, rotate=90, anchor=south] at (axis description cs:-0.085,0.18) {less robust};
        \draw[robarrow] (axis description cs:-0.085,0.66) -- (axis description cs:-0.085,0.98);
        \node[roblabel, rotate=90, anchor=south] at (axis description cs:-0.085,0.82) {more robust};
    \nextgroupplot[xmin=0, xmax=0.7, ymin=0, ymax=0.7]
        \addplot[fill=green!30, draw=none]
            coordinates {(0,0) (0.7,0) (0.7,0.7) (0,0.7)} \closedcycle;
        \addplot[fill=gray!18, draw=none] coordinates {(0.0000,0.0000) (0.0000,0.7000) (0.0223,0.7000) (0.0374,0.4180) (0.0524,0.2980) (0.0675,0.2315) (0.0826,0.1893) (0.0976,0.1601) (0.1127,0.1387) (0.1277,0.1223) (0.1428,0.1094) (0.1579,0.0990) (0.1729,0.0904) (0.1880,0.0831) (0.2030,0.0770) (0.2181,0.0716) (0.2332,0.0670) (0.2482,0.0629) (0.2633,0.0593) (0.2783,0.0561) (0.2934,0.0533) (0.3085,0.0507) (0.3235,0.0483) (0.3386,0.0461) (0.3536,0.0442) (0.3687,0.0424) (0.3837,0.0407) (0.3988,0.0392) (0.4139,0.0378) (0.4289,0.0364) (0.4440,0.0352) (0.4590,0.0340) (0.4741,0.0330) (0.4892,0.0319) (0.5042,0.0310) (0.5193,0.0301) (0.5343,0.0292) (0.5494,0.0284) (0.5645,0.0277) (0.5795,0.0270) (0.5946,0.0263) (0.6096,0.0256) (0.6247,0.0250) (0.6398,0.0244) (0.6548,0.0239) (0.6699,0.0233) (0.6849,0.0228) (0.7000,0.0223) (0.7000,0.0000)} \closedcycle;
        \addplot[green!55!black, very thick, domain=0.0223:0.7, samples=200]
            {1/(64*x)};
        \node[circle, fill=black, inner sep=1.1pt] at (axis cs:0,0) {};
        \node[anchor=south west, font=\small, inner sep=2.5pt,
              fill=white, fill opacity=0.75, text opacity=1]
              at (axis cs:0.04,0.04) {Nash};
        \draw[robarrow] (axis description cs:0.34,-0.085) -- (axis description cs:0.02,-0.085);
        \node[roblabel, anchor=north] at (axis description cs:0.18,-0.085) {less robust};
        \draw[robarrow] (axis description cs:0.66,-0.085) -- (axis description cs:0.98,-0.085);
        \node[roblabel, anchor=north] at (axis description cs:0.82,-0.085) {more robust};
    \end{groupplot}
        \path (group c1r1.north west) -- (group c2r1.north east)
        node[midway, above=3pt, font=\footnotesize, inner sep=0pt]
        {{\setlength{\fboxsep}{0pt}%
          \fcolorbox{green!55!black}{green!30}{\rule{0pt}{1.4ex}\rule{2.4ex}{0pt}}~unique \& tractable%
          \hspace{1.2em}%
          \fcolorbox{black!45}{gray!18}{\rule{0pt}{1.4ex}\rule{2.4ex}{0pt}}~no certificate%
          \hspace{1.2em}%
          \fcolorbox{red!45}{red!15}{\rule{0pt}{1.4ex}\rule{2.4ex}{0pt}}~ill-posed}};
    \end{tikzpicture}
    \vspace{-10.0pt}
    \caption{Certified region (green) of $\chi^1\chi^2>1$ in the
    $(1/m^1,1/m^2)$ plane, for two-player linear-quadratic data: running costs
    $\norm{x_t}^2_{Q}+\norm{u_t^i}^2_{R^i}
    -\norm{d_t^i}^2_{S^i}$ and penalty
    $\norm{d_t^i-u_t^{-i}}^2_{M^i}$, with
    $Q=I$, $R^i=r^iI$, $S^i=sI$, $M^i=m^iI$, $r^1=r^2=1$, and
    $\sigma_{\min}(G^i)=\sigma_{\max}(G^{-i})=1$.
    \emph{Left} ($s=0$): a bounded island. The certificate fails near Nash,
    where the interference swamps the curvatures, and past the solvability
    limit (red), where the adversary's problem loses concavity.
    \emph{Right} ($s=1$): the direct disturbance penalty enters only through
    $\mu_D^i$, raising it. Here $\mu_D^i=2m^i>0$ for every $m^i>0$: the
    second failure mode disappears, the region is unbounded, and arbitrarily
    strong robustness is certified, approaching robust control.}
    \label{fig:reweight}
    \vspace{-2.5ex}
\end{figure}

This two-player certificate is optimal among positive reweightings based on the stated bounds; Figure~\ref{fig:reweight} shows the certified region. More generally, the two-player condition is $(c^{12}c^{21})^2 < 16\mu_U^1\mu_D^1\mu_U^2\mu_D^2$.

For $N>2$, the reweighting~\eqref{eq:weights} need not be optimal, but optimal weights can be obtained by solving a linear matrix inequality in~$\lambda$, which we omit for brevity. Finally, if $c_{\mathrm{out}}^ic_{\mathrm{in}}^j<4\mu_D^i\mu_U^j$ for all $i\neq j$, equal weights also work, so $F$ itself is strongly monotone. The certificate can also be met by design: the robustness parameters directly control $c^{ij}$ and $\overline\ell^i$, and, whenever $\mu_U^k>0$, a sufficiently strong direct disturbance penalty combined with a sufficiently small deviation penalty satisfies~\eqref{eq:monotonicity_condition}.

\section{The closed-loop linear-quadratic setting}\label{sec:closed_loop}
Throughout this section, we consider the (unconstrained) linear-quadratic setting, i.e., Assumption~\ref{ass:cl} is in force.
Players now choose
state-feedback laws: the stage policy set $\Pi_t^i$ is the set of measurable
functions from the state history $(x_0,\ldots,x_t)\in(\reals^n)^{t+1}$ to
the stage input space $\reals^{n_i}$, and the controller is $u_t^i = \pi_t^i(x_0, \ldots, x_t)$.
We write $\Pi^i=\Pi_0^i\times\cdots\times\Pi_{T-1}^i$ for the policy set of
player $i$ and $\Sigma^i=\prod_{j\neq i}\Pi^j$ for the joint policy set available
to their fictitious adversary.
The nominal reference actions in player $i$'s
deviation penalty are evaluated along the player-specific fictitious trajectory
of Section~\ref{sec:setting}, i.e.,
$u_t^{-i}=\pi_t^{-i}(x_0^{(i)},\ldots,x_t^{(i)})$.

\begin{definition}[closed-loop strategically robust dynamic equilibrium]\label{def:sre:cl}
    A tuple of policies $(\bar\pi^{1},\ldots,\bar\pi^N)\in\Pi^1\times\cdots\times\Pi^N$ forms a closed-loop strategically robust dynamic equilibrium if for all players $i$
    \begin{equation}\label{eq:sre:cl}
        J^i(\bar\pi^i,\bar\pi^{-i})
        \leq
        J^i(\pi^i,\bar\pi^{-i})
        \qquad
        \forall \pi^i\in\Pi^i.
    \end{equation}
\end{definition}

Finally, we call a strategically robust dynamic equilibrium \emph{subgame perfect} if the equilibrium condition holds in every continuation game after every admissible state history; existence, uniqueness, and computation are studied in this sense.

Definitions~\ref{def:sre:ol} and~\ref{def:sre:cl} impose the same
inequality on the strategically robust cost~\eqref{eq:robust_cost}: no
unilateral deviation of player $i$ lowers their robust cost. They differ only in
the set over which that deviation ranges---input sequences in $\mathcal U^i$, or
state-feedback policies in $\Pi^i$---and, correspondingly, in the set available
to the fictitious adversary.

The augmented cost of player $i$ therefore reads
\begin{multline*}
 J_\mathrm{aug}^i(\pi^i,\sigma^i,\pi^{-i}) = \norm{x_T}^2_{Q^i_T} +
    \\
    \sum_{t=0}^{T-1} \norm{x_t}^2_{Q_t^i} + \norm{u_t^i}^2_{R_t^i}
    - \norm{d_t^i}^2_{S_t^i}
    - \norm{d_t^{i}-u_t^{-i}}^2_{M_t^i}
\end{multline*}
subject to the system dynamics, with the adversary maximizing the same expression. The term $\norm{u_t^{-i}}^2_{S_t^i}$ is constant in player $i$'s open-loop decision and so does not affect its open-loop best response; under feedback it
acts indirectly through the state trajectory. We keep it as a direct regularizer
of the fictitious adversary, and setting $S_t^i=0$ recovers the standard LQ
running cost; Section~\ref{sec:hinf} explains its role in the connection with
$\mathcal{H}_\infty$ control. Similarly, as in Section~\ref{sec:ol_static}, we can reformulate the problem as a simultaneous game using the augmented cost.

\begin{definition}[equilibrium of the closed-loop augmented game]\label{def:aug:cl}
    The policies $(\bar\pi^1,\bar\sigma^{1}, \ldots, \bar\pi^N,\bar\sigma^{N})$ are an equilibrium of the augmented game if for all players $i$ we have
    \begin{equation}\label{eq:aug:cl}
        J_\mathrm{aug}^i(\bar \pi^{i}, \sigma^i, \bar\pi^{-i})
        \leq
        J_\mathrm{aug}^i(\bar \pi^{i}, \bar\sigma^i, \bar\pi^{-i})
        \leq
        J_\mathrm{aug}^i(\pi^i, \bar\sigma^i, \bar\pi^{-i})
    \end{equation}
    for all policies $\pi^i\in\Pi^i$ and $\sigma^i\in\Sigma^i$.
\end{definition}

\begin{lemma}[equivalence in the closed-loop setting]\label{lemma:equivalence:cl}\,
\begin{enumerate}
    \item [\textnormal{(a)}]
    \emph{Necessity.}
    If $(\bar\pi^1,\bar\sigma^{1},\ldots,\bar\pi^N,\bar\sigma^{N})$ is an
    equilibrium of the augmented game, then $(\bar\pi^1,\ldots,\bar\pi^N)$ is a
    closed-loop strategically robust dynamic equilibrium and each
    $\bar\sigma^i$ is a worst-case policy in~\eqref{eq:robust_cost}, i.e., a
    best response to $\bar\pi^i$.
    \item [\textnormal{(b)}]
    \emph{Sufficiency.}
    Suppose
    $(\bar\pi^1,\ldots,\bar\pi^N)$ is a closed-loop strategically robust
    dynamic equilibrium and that, for each player $i$, the zero-sum game
    $J_{\mathrm{aug}}^i(\cdot,\cdot,\bar\pi^{-i})$ admits a saddle point.
    Then, there exist worst-case policies $\bar\sigma^1,\ldots,\bar\sigma^N$
    such that
    $(\bar\pi^1,\bar\sigma^{1},\ldots,\bar\pi^N,\bar\sigma^{N})$ is an
    equilibrium of the augmented game.
\end{enumerate}
\end{lemma}

Lemma~\ref{lemma:equivalence:cl} is proved in Section~\ref{sec:proofs_cl} by specializing the general equivalence principle of Lemma~\ref{lemma:equivalence:gen} to the closed-loop strategy sets.

\subsection{Relating Nash equilibria and $\mathcal{H}_\infty$ control}\label{sec:hinf}
Before the equilibrium characterization, we explain how strategic robustness
connects two classical extremes: the Nash equilibrium (in particular, the best-response problem with the other players' policies as given) and $\mathcal H_\infty$ control (which
treats their inputs as adversarial disturbances). This connection will also drive the infinite-horizon analysis in
Section~\ref{sec:infinite}.

Fix a player $i$, take $S_t^i=s^iI$ with $s^i>0$, and first set
$M_t^i\equiv0$. The opponents' policies then drop out
of~\eqref{eq:robust_cost}, and player $i$ is left with
\begin{equation}\label{eq:hinf:problem}
    \inf_{\pi^i}\ \sup_{\sigma^i}\ \norm{x_T}^2_{Q_T^i}
    +\sum_{t=0}^{T-1}\norm{x_t}^2_{Q_t^i}+\norm{u_t^i}^2_{R_t^i}
    -s^i\norm{d_t^i}^2,
\end{equation}
subject to $x_{t+1}=A_tx_t+B_t^iu_t^i+B_t^{-i}d_t^i$. This is the setting of the classical
discrete-time $\mathcal H_\infty$ problem: $d^i$ is an unrestricted disturbance
entering through $B_t^{-i}$, and $s^i$ is the level at which its energy is
charged, so the smaller $s^i$, the more powerful the disturbance and the more
conservative the controller. The problem has a finite value, attained by a
linear state feedback, only when $s^i$ is large enough; let $s_{\mathrm{att}}^i$ denote
the infimum of the levels for which the value is finite.
Now reinstate the penalty $m^i\norm{d_t^i-u_t^{-i}}^2$ and complete the square
to get
\begin{align}
s^i\norm{d_t^i}^2+m^i\norm{d_t^i-u_t^{-i}}^2
&=(s^i+m^i)\left\|d_t^i-\frac{m^i}{s^i+m^i}u_t^{-i}\right\|^2
\nonumber\\
&\quad+\frac{s^im^i}{s^i+m^i}\norm{u_t^{-i}}^2.
\label{eq:complete_square}
\end{align}

The adversary faces a single
magnitude penalty with weight $s^i+m^i$, centered at the anchor
$\frac{m^i}{s^i+m^i}u_t^{-i}$.
For fixed $s^i>0$, the limits $m^i\to0$ and $m^i\to\infty$ yield,
respectively, the setting of the classical $\mathcal H_\infty$ problem and the nominal best-response problem associated with $S_t^i=s^iI$. Thus, the standard Nash
problem ($s^i=0$, $m^i\to\infty$) and the $\mathcal H_\infty$ problem
($s^i>s_{\mathrm{att}}^i$, $m^i\to0$) are distinct limits of the same two-parameter
family.

\subsection{The finite-horizon case}\label{sec:finite}
Over a finite horizon the equilibrium is characterized by a coupled backward recursion, along which two kinds of condition appear. This mirrors the standard theory of dynamic linear-quadratic games.

\begin{theorem}[existence, uniqueness, and computation in the closed-loop setting]\label{thm:sre}
Starting from $P_T^i=Q_T^i$, define, for each player $i$ and time $t$,
\begin{align*}
    H_t^{u^i,u^i}
    &\coloneqq R_t^i+(B_t^i)^\top P_{t+1}^iB_t^i,\\
    H_t^{u^i,d^i}
    &\coloneqq(B_t^i)^\top P_{t+1}^iB_t^{-i},\\
    H_t^{d^i,u^i}
    &\coloneqq-(B_t^{-i})^\top P_{t+1}^iB_t^i,\\
    H_t^{d^i,d^i}
    &\coloneqq M_t^i+S_t^i-(B_t^{-i})^\top P_{t+1}^iB_t^{-i},\\
    G_t^{u^i}
    &\coloneqq-(B_t^i)^\top P_{t+1}^iA_t,\qquad
    G_t^{d^i}\coloneqq(B_t^{-i})^\top P_{t+1}^iA_t.
\end{align*}
Let $H_t$ be the block matrix indexed by
$(u^1,\ldots,u^N,d^1,\ldots,d^N)$ whose nonzero blocks are
\begin{equation*}
\begin{aligned}
    (H_t)_{u^i,u^i}&=H_t^{u^i,u^i},&
    (H_t)_{u^i,d^i}&=H_t^{u^i,d^i},\\
    (H_t)_{d^i,u^i}&=H_t^{d^i,u^i},&
    (H_t)_{d^i,d^i}&=H_t^{d^i,d^i},\\
    (H_t)_{d^i,u^j}&=-M_t^{i,j},&&j\neq i.
\end{aligned}
\end{equation*}
Here, $M_t^{i,j}$ denotes the block columns of $M_t^i$ associated with
player $j$ in the stacked vector $u_t^{-i}$; all remaining blocks are zero.
Let
\begin{align*}
K_t&\coloneqq\operatorname{col}(K_t^1,\ldots,K_t^N),&
L_t&\coloneqq\operatorname{col}(L_t^1,\ldots,L_t^N),\\
G_t^u&\coloneqq\operatorname{col}(G_t^{u^1},\ldots,G_t^{u^N}),&
G_t^d&\coloneqq\operatorname{col}(G_t^{d^1},\ldots,G_t^{d^N}),
\end{align*}
and define the feedback gains by
\begin{equation}\label{eq:sre_policies}
    \operatorname{col}(K_t,L_t)
    =H_t^{-1}\operatorname{col}(G_t^u,G_t^d).
\end{equation}
For
${A^i_{\mathrm{CL}}}_t\coloneqq A_t+B_t^iK_t^i+B_t^{-i}L_t^i$, update
\begin{equation}\label{eq:sre_riccati}
\begin{aligned}
P_t^i
&=Q_t^i+(K_t^i)^\top R_t^iK_t^i-(L_t^i)^\top S_t^iL_t^i\\
&\quad-(L_t^i-K_t^{-i})^\top M_t^i(L_t^i-K_t^{-i})\\
&\quad+({A^i_{\mathrm{CL}}}_t)^\top P_{t+1}^i{A^i_{\mathrm{CL}}}_t.
\end{aligned}
\end{equation}
Suppose that, along this backward recursion,
$H_t^{u^i,u^i}\succ0$, $H_t^{d^i,d^i}\succ0$ for every $i$, and
$H_t$ is invertible at every time $t$. Then the policies
$\bar\pi_t^i(x)=K_t^ix$ form the unique subgame-perfect strategically robust
dynamic equilibrium. Moreover, at every continuation history the unique
worst-case action in~\eqref{eq:robust_cost} is $d_t^i=L_t^ix_t$.
\end{theorem}

Theorem~\ref{thm:sre}, proved in Section~\ref{sec:proofs_cl}, shows
that the strategically robust equilibrium retains the familiar structure of
finite-horizon LQ Nash equilibria: the policies are linear and Markovian, and
their gains are obtained through a coupled backward Riccati recursion. The
additional adversarial gains quantify the policy deviations against which each
player is protected.
The stagewise conditions $H_t^{u^i,u^i}\succ0$ and $H_t^{d^i,d^i}\succ0$ keep each stage problem well posed, and they are close to necessary: a negative eigenvalue in either lets the corresponding player or adversary drive the stage cost-to-go to $\mp\infty$. The remaining condition, invertibility of $H_t$, is the only one coupling the players; we return to it in Proposition~\ref{prop:invertibility}.
The method of the proof---backward dynamic programming on a stagewise saddle---is not
tied to the linear-quadratic case; what is tied to it is the closed form. Only here
do the value functions stay quadratic, and each stage reduces to a linear system,
so that the whole construction collapses to a Riccati recursion.

\begin{example}[motivating, continued]
    Consider the setup of Example~\ref{example:motivating} extended to the framework of strategically robust agents. We set the robustness levels to $m^1_t = 1.2$ and $m^2_t \to +\infty$, i.e., player 1 seeks protection against misspecified policies of player 2.
    Under the resulting strategically robust dynamic equilibrium, for the majority of the horizon, the state evolves according to $x_{t + 1} \approx 0.464 x_t$, remaining stable even under perturbations in the feedback policy of player 2.
\end{example}

Invertibility of $H_t$ ties the players' parameters together; such conditions also appear in standard linear-quadratic games, e.g.,~\cite{bacsar1998dynamic}.
In our case, the reweighting certificate of Theorem~\ref{thm:monotonicity} also gives a sufficient condition at each stage: self-curvatures of players and adversaries must dominate the strategic interference induced by the penalties.

\begin{proposition}[invertibility of $H_t$]\label{prop:invertibility}
    At time $t$, define, for each player $i$,
    \begin{align*}
        \mu_U^i &\coloneqq
        \lambda_{\min}(R_t^i + (B_t^i)^\top P_{t+1}^i B_t^i), \\
        \mu_D^i &\coloneqq
        \lambda_{\min}(S_t^i + M_t^i - (B_t^{-i})^\top P_{t+1}^i B_t^{-i}),
    \end{align*}
    and let $c^{ij}\coloneqq\sigma_{\max}(M_t^{i,j})$ for $j\neq i$.
    Form the ratios in~\eqref{eq:cycle_constants} with these quantities.
    If $\mu_U^i,\mu_D^i>0$ for every $i$ and
    \eqref{eq:monotonicity_condition} holds, then $H_t$ is invertible.
    In particular, if $M_t^i=mI\succ0$ for every $i$, then
    \begin{equation*}
        \chi^i=\frac{4\mu_U^i\mu_D^i}{m^2(N-1)^2},
    \end{equation*}
    and the condition is applied to the two smallest of these ratios.
\end{proposition}

This plays the role here that Theorem~\ref{thm:monotonicity} plays in the open loop: player and adversary curvatures at stage $t$ must dominate the interference the deviation penalties create between them, and since $P_{t+1}^i$ is known when stage $t$ is processed, the condition is checkable as the recursion runs.

\subsection{Towards the infinite horizon}\label{sec:infinite}

We now turn our attention to the infinite-horizon case.
Outside the zero-sum case, infinite-horizon linear-quadratic games of this kind remain largely
open. A stationary equilibrium is characterized by a system of algebraic
Riccati equations, nonlinear in the unknowns and coupled across all players.
Such a system need not admit a solution, and when it does the solution need not
be unique~\cite{engwerda2005lq,salizzoni2025nash}.
Nor is there an algorithm known to reach one from an
arbitrary start; e.g., policy-gradient schemes can fail to converge even locally in
$N$-player linear-quadratic games~\cite{mazumdar2020policy}.
We now show that strategic robustness can alleviate these issues, provided that players are sufficiently robust.

As is standard over an infinite horizon, we take the data time-invariant. The state evolves as $x_{t+1}=Ax_t+\sum_{j=1}^N B^j u_t^j$, and player $i$ minimizes the strategically robust cost
\begin{equation*}
J^i(\pi^i, \pi^{-i})=\sup_{\sigma^i}\sum_{t=0}^{\infty}\norm{x_t}^2_{Q^i}+\norm{u_t^i}^2_{R^i} -\norm{d_t^i}^2_{S^i}-\norm{d_t^i-u_t^{-i}}^2_{M^i},
\end{equation*}
subject to the fictitious dynamics
$x_{t+1}=Ax_t+B^iu_t^i+B^{-i}d_t^i$ generated by the policies $(\pi^i, \sigma^i)$, where, as above,
$d_t^i\in\reals^{\sum_{j\neq i}n_j}$ is player $i$'s adversary acting through
the stacked opponent channel $B^{-i}$. The supremum is over all causal measurable adversarial policies, and infinite-horizon sums are interpreted as the limsup of their finite partial sums. We take $S^i=s^iI\succeq0$ and
$M^i=m^iI\succ0$. We now define the equilibrium concept:

\begin{definition}[stationary linear strategically robust dynamic equilibrium]\label{def:sre:inf}
    A tuple of policies $(\bar\pi^1,\ldots,\bar\pi^N)$ is a stationary linear
    strategically robust dynamic equilibrium if there exist gains
    $K^1,\ldots,K^N$ with $\bar\pi_t^i(x)=K^ix$ for every player $i$ and every
    time $t$ such that
    \begin{equation}\label{eq:sre:inf}
        J^i(\bar\pi^i,\bar\pi^{-i})
        \leq
        J^i(\pi^i,\bar\pi^{-i})
    \end{equation}
    for every initial state $x_0$, for every $i$, and every $\pi^i=(\pi_0^i,\pi_1^i,\ldots)$ whose stage maps
    $\pi_t^i$ are measurable functions of the state history.
\end{definition}

Here, $K^i$ is stationary (time-invariant), while the deviations $\pi^i$ are unrestricted.
We begin with a verification statement: if the gains satisfy the following assumptions, we have a strategically robust dynamic equilibrium.

\begin{proposition}[verification in the infinite horizon]\label{prop:inf_verify}
Suppose there exist $P^1,\ldots,P^N\succeq0$ and gains $K,L$. For each player
$i$, set
\begin{equation*}
\mathcal H^i(P^i)\coloneqq
\begin{bmatrix}
H^{u^i,u^i} & H^{u^i,d^i}\\
H^{d^i,u^i} & H^{d^i,d^i}
\end{bmatrix},
\end{equation*}
where the blocks are the same as in Theorem~\ref{thm:sre} with a time-independent $P^i$. Suppose that:
\begin{enumerate}
\item[\textnormal{(i)}] for every $i$, the symmetric part of
$\mathcal H^i(P^i)$ is positive definite and
\begin{equation*}
\mathcal H^i(P^i)
\begin{bmatrix}K^i\\ L^i\end{bmatrix}
=
\begin{bmatrix}
-(B^i)^\top P^iA\\
(B^{-i})^\top P^iA+m^iK^{-i}
\end{bmatrix};
\end{equation*}
\item[\textnormal{(ii)}] the gains satisfy
\begin{equation}\label{eq:inf:verify_domain}
Q^i-\frac{s^im^i}{s^i+m^i}(K^{-i})^\top K^{-i}\succ0
\end{equation}
and render every $A^i_{\mathrm{CL}}\coloneqq A+B^iK^i+B^{-i}L^i$ Schur;
\item[\textnormal{(iii)}] each $P^i$ satisfies
\begin{equation}\label{eq:inf:care}
\begin{aligned}
P^i&=Q^i+(K^i)^\top R^iK^i-(L^i)^\top S^iL^i\\
&-(L^i-K^{-i})^\top M^i(L^i-K^{-i})
+(A^i_{\mathrm{CL}})^\top P^iA^i_{\mathrm{CL}}.
\end{aligned}
\end{equation}
\end{enumerate}
Then $(\bar \pi^1,\ldots,\bar \pi^N)$ with $\bar \pi_t^i(x) = K^i x$ is a stationary linear strategically robust dynamic
equilibrium, and $\sigma^i_t(x) = L^ix$ is the policy that achieves the inner supremum.
\end{proposition}

The bound~\eqref{eq:inf:verify_domain} makes the centered state weights positive, while conditions~(i)--(iii) are the standard stabilizing $\mathcal H_\infty$ verification conditions~\cite[Ch.~3]{bacsar2008h}. The next theorem constructs such solutions and proves contraction on an invariant region.
The only structural input is that each player's problem be solvable in the fully decoupled case when $m^i = 0$. This case is equivalent to $N$ independent $\mathcal{H}_{\infty}$ control problems.

\begin{assumption}[per-player feasibility]\label{ass:feas}
For each player $i$, it holds that $Q^i\succ0$, $R^i\succ0$, and
$S^i=s^iI\succ0$. The decoupled problem obtained by setting $m^i=0$ admits the
usual infinite-horizon $\mathcal H_\infty$ stabilizing solution, with value
matrix $P_\infty^i \succ 0$ and stationary optimal-controller and worst-response gains $K_\infty^i,L_\infty^i$.
Explicitly,
\begin{equation}\label{eq:inf:nominal_stabilizing}
 s^iI-(B^{-i})^\top P_\infty^iB^{-i}\succ0,
 F_\infty^i\coloneqq A+B^iK_\infty^i+B^{-i}L_\infty^i
\end{equation}
is Schur. Equivalently, $s^i$ lies in the strict solvability region of the decoupled
state-feedback $\mathcal H_\infty$ problem (in particular, strictly above its
optimal attenuation threshold).
\end{assumption}

Denote by $\mathrm{BR}^i(K^{-i})$ the best response of player $i$ when the others play $K^{-i}$. We first show that this best response exists and is unique on a ball of gains around the origin, and then that the resulting map $\mathrm{BR}$ is a contraction there. The equilibrium is not merely shown to exist, but is also computed by iterating the best response $\mathrm{BR}$.

\begin{theorem}[existence, uniqueness, and computation in infinite horizon]\label{thm:inf}
Let Assumption~\ref{ass:feas} hold. Fix $\bar k_0>0$ independently of
$m^1,\ldots,m^N$ and let
$\mathcal B_0\coloneqq\{K:\max_j\norm{K^j}\le\bar k_0\}$.
\begin{enumerate}
\item[\textnormal{(a)}] For all sufficiently small $m^1,\ldots,m^N>0$, for every
$K\in\mathcal B_0$ and every $i$, the best
response $\mathrm{BR}^i(K^{-i})$ exists, is unique, and is continuously
differentiable, so that
\begin{equation}\label{eq:inf:C}
K_0^i\coloneqq\mathrm{BR}^i(0),
\quad
C_{ij}\coloneqq\!\!\sup_{K\in\mathcal B_0}\!\!
\norm{D_{K^j}\mathrm{BR}^i(K^{-i})}<\infty
\end{equation}
for $j\neq i$; set $C_{ii}\coloneqq0$ and
$C\coloneqq[C_{ij}]\in\reals_{\ge0}^{N\times N}$.

\item[\textnormal{(b)}] If
\begin{equation}\label{eq:inf:conditions}
\rho(C)<1
\quad\text{and}\quad
\norm{K_0^i}+\bar k_0\textstyle\sum_{j\neq i}C_{ij}\le\bar k_0
\ \ \forall i,
\end{equation}
then $\mathrm{BR}(\mathcal B_0)\subseteq\mathcal B_0$ and $\mathrm{BR}$ has a
unique fixed point $K_\star\in\mathcal B_0$, which is the unique stationary
linear strategically robust dynamic equilibrium in $\mathcal B_0$. For every
$K^{(0)}\in\mathcal B_0$, $\mathrm{BR}^{(n)}(K^{(0)})\to K_\star$ geometrically,
with asymptotic root rate at most $\rho(C)$.

\item[\textnormal{(c)}] Fix $\bar k_0>\max_i\norm{K_\infty^i}$ independently of
$m^1,\ldots,m^N$. As $\max_im^i\to0$,
\begin{equation}\label{eq:inf:small_m_rates}
K_0^i=K_\infty^i+O(m^i),
\quad
C_{ij}=O(m^i),
\end{equation}
with constants independent of $m^1,\ldots,m^N$. Hence~\eqref{eq:inf:conditions} holds, and~(a) and~(b)
apply, for all sufficiently small $m^1,\ldots,m^N$.
\end{enumerate}
\end{theorem}

In words, Part~(a) is a local perturbation result from
the decoupled stabilizing $\mathcal H_\infty$ solutions. The second condition
in~\eqref{eq:inf:conditions} keeps the best responses inside $\mathcal B_0$,
while $\rho(C)<1$ makes the iteration converge and bounds its asymptotic root
rate. Entry $C_{ij}$ measures how far player $i$'s optimal gain moves when
opponent $j$'s gain moves, so $C$ is the network of cross-player sensitivities,
and $\rho(C)<1$ asks that these sensitivities not reinforce one another around
any cycle of players rather than that each be small. Part~(c) delivers the
invariance and contraction conditions at small $m^i$, since the
deviation penalty is the only channel coupling the players' fictitious
problems.

\begin{figure*}[!t]
    \centering
    \begin{minipage}[t]{0.8\textwidth}
        \vspace{0pt}
        \centering
        \includegraphics[scale=0.95, trim=11bp 4bp 216bp 14bp, clip]{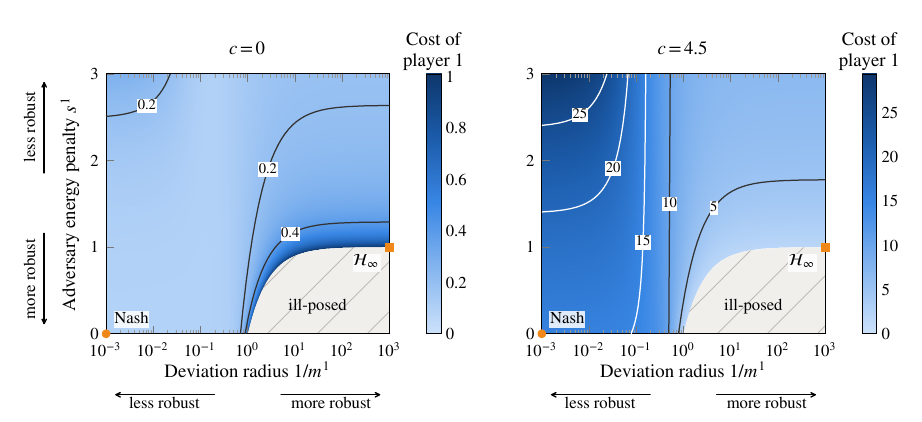}%
        \hspace{14pt}%
        \includegraphics[scale=0.95, trim=252bp 4bp 9bp 14bp, clip]{figures/robustness_costs.pdf}
        \par\vspace{-6pt}
        \caption{Realized physical cost of player 1 under perturbation budgets $c \in \{0, 4.5\}$. The most robust policies incur greater zero-perturbation cost but degrade substantially less as the perturbation budget increases. Furthermore, there exist points that have relatively low costs under both regimes; for example $(s^1, m^1) = (s_{\mathrm{att}}, s_{\mathrm{att}}) = (1, 1)$.}
        \label{fig:adversarial_costs}
    \end{minipage}%
    \hfill
    \begin{minipage}[t]{0.19\textwidth}
        \vspace{0pt}
        \centering
        \setlength{\tabcolsep}{4pt}%
        \captionof{table}{Realized physical-cost statistics comparing the nominal and strategically robust policies with $m^1=0.1$}
        \label{tab:cost_statistics}
        \footnotesize
        \begin{tabular}{ccc}
            \toprule
            Percentile & Nom. & SR \\
            \midrule
            5\%  & 0.07  & 0.65  \\
            25\% & 1.18  & 1.45  \\
            50\% & 5.18  & 4.40  \\
            75\% & 15.09 & 11.67 \\
            95\% & 43.74 & 33.03 \\
            \bottomrule
        \end{tabular}
    \end{minipage}
    \vspace{-15pt}
\end{figure*}

\section{Numerical experiments}\label{sec:experiments}
We begin with the closed-loop experiments and turn afterwards to the open-loop experiments.

\subsection{Closed-loop experiments}

In the closed-loop setting, the equilibrium is obtained directly from the coupled backward Riccati recursion of Theorem~\ref{thm:sre}. At each stage, we numerically verify the positive-definiteness conditions on $H_t^{u^i,u^i}$ and $H_t^{d^i,d^i}$, as well as the invertibility of the joint matrix $H_t$.

\subsubsection{One-dimensional collaborative game}

Two agents collaborate to steer a scalar system from $x_0=1$ to the target $x_T=0$ over $T=3$ stages, with integrator dynamics $x_{t+1} = x_t + u_t^1 + u_t^2$ and control cost $x_3^2 + \sum_{t=0}^{T-1} (u_t^i)^2$.
Across all settings, the robustness matrices are time-invariant
$M_t^i=m^iI, S_t^i = s^i I$. When only player~1 is robust, we set the other player to be nominal at $s^2=0, m^2 \to \infty$, and specify $(s^1, m^1)$ for player 1. We refer to $s_{\mathrm{att}}$ as the infimum of $s^i$ such that the individual problem is solvable at $m^i = 0$. When both players have equal robustness, we choose $s^1=s^2=s_{\mathrm{att}}$, and set $m = m^i >0$.

\paragraph*{Robustness to adversarial perturbations}
We make only player~1 robust and assess the performance of strategically robust policies prescribed to player~1 under the entire $(s^1, m^1)$ spectrum. Player~2 uses its Nash equilibrium strategy, which we perturb toward player~1's adversarial policy:
\begin{equation*}
    \pi^{2,\mathrm{adv}}_t(x_t)
    =
    \pi^{2,\mathrm{NE}}_t(x_t) - \mathrm{clip}(\pi^{2,\mathrm{NE}}_t(x_t) - d^{1,2}_{\mathrm{robust}, t}(x_t), -c,c),
\end{equation*}
where $d^{1,2}_{\mathrm{robust}, t}$ is the adversarial policy obtained when synthesizing player~1's robust controller recovered from $m^1 = 0.1, s^1 = s_{\mathrm{att}}$. The parameter $c\geq0$ is the perturbation budget. Figure~\ref{fig:adversarial_costs} displays the full spectrum of parameterizations. For a low value of $c$, those near the solvability boundary pay a higher cost, but most are distributed evenly. When the perturbation budget becomes large, the strategies near the boundary, which are also the most robust ones, perform better.

\paragraph*{Robustness to random perturbations}
We only make player~1 robust, and we corrupt player~2's feedback by a constant random bias,
$\pi_t^{2,\mathrm{bias}}(\cdot)=\pi_t^{2,\mathrm{Nom}}(\cdot)+b$, where $\pi_t^{2,\mathrm{Nom}}(\cdot)$ is the nominal equilibrium policy attained from $s^1=s_{\mathrm{att}}, m^1 \to \infty$, and
$b\sim\mathcal N(0,4)$.
Through Monte Carlo simulations with $10^5$ samples, we compare in Table~\ref{tab:cost_statistics} the nominal policy of player~1 and the robust policy with $m^1=0.1, s^1 = s_{\mathrm{att}}$.
Relative to the nominal policy, the robust policy reduces the 95th and 75th percentiles by 24\% and 23\% and the median by 15\%, while raising the lower percentiles slightly.

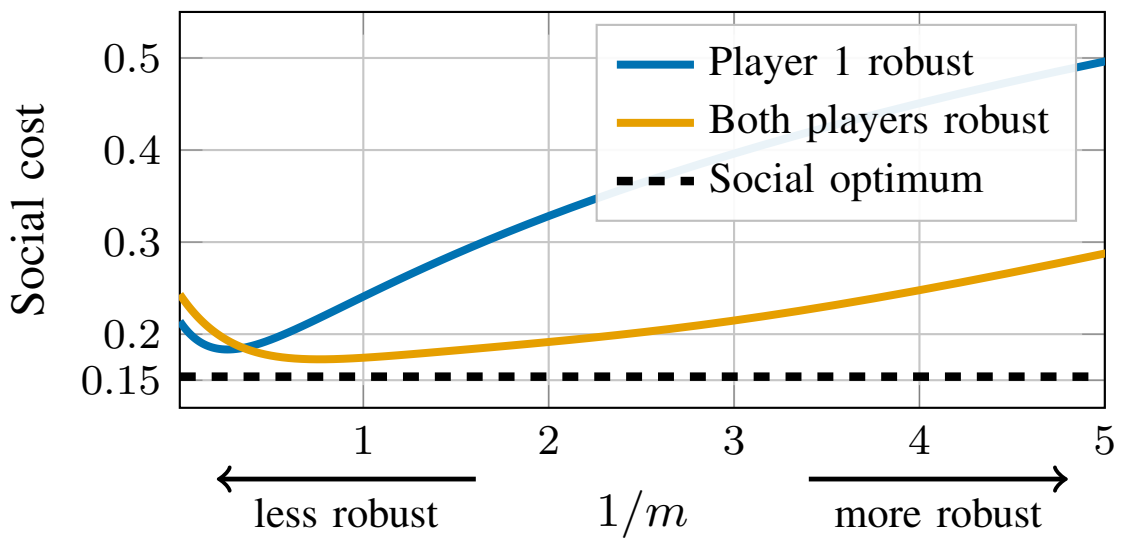
\begin{figure}[!t]
    \centering
    \scalebox{\figscale}{\input{figures/free_lunch_plot}}
    \vspace{-1.0\baselineskip}
    \caption{Social cost, defined as the sum of the agents' realized physical costs, as a function of the strategic-robustness level $1/m$ when only player 1 is robust (blue) and when both players are (yellow).}
    \label{fig:collaborative:freelunch}
    \vspace{-2.5ex}
\end{figure}

\paragraph*{Robustness induces collaboration}
Our last experiment for this game concerns the social cost, the sum of all agents' costs.
We compare the case where only player~1 is robust ($s^1 = s_{\mathrm{att}}, m^1 = m$) to the case when both are; the results appear in Figure~\ref{fig:collaborative:freelunch}.
The social cost decreases over part of the robustness range---a ``free-lunch'' effect related to the coordination-via-robustification phenomenon observed in static games~\cite{lanzetti2025strategically} and used in~\cite{qu2026training} for training collaborative AI agents. Unlike the usual nominal-performance tradeoff in single-agent robust optimization, the strategic response of the other player can improve nominal aggregate performance. In this example, robust players protect against the possibility that the other agent exerts less control effort by exerting more themselves, which steers the system closer to the target and lowers the social cost.

\subsubsection{Network consensus and control}

Our second closed-loop experiment is a multi-agent control problem on a graph~\cite{pasqualetti2014controllability}: a star network of $N=5$ nodes (Figure~\ref{fig:network:star}), with node~1 at the center and nodes~2--5 as leaves.
The system state of node $i$ evolves according to the dynamics
$
    x_{t+1}^i
    =
    \frac{2}{3} x_t^i
    +
    \frac{1}{3\abs{\mathcal N^i}} \sum_{j\in \mathcal N^i}x_t^j
    +
    u_t^i,
$
where $\mathcal N^i$ is the set of neighbors of node $i$ and the initial condition is $x_0^i=0$.
Each agent $i$ aims to minimize control effort while reaching the target state $x_\mathrm{target}=5$, so that
$
    J_{\mathrm{nom}}^i(\pi^i,\pi^{-i}) = \sum_{t=0}^{T-1} (u_t^i)^2 + (x^i_{T} - 5)^2,
$
over a horizon of $T=20$. The affine target can be represented in our formulation through the standard augmentation of the state by a constant component. The robustness matrices are isotropic and time-invariant across players. If we are making only the center robust, we set $s^{\mathrm{center}} = s_{\mathrm{att}}^{\mathrm{center}}, s^{\mathrm{leaf } j} = 0, m^{\mathrm{leaf }j} \to \infty$, and set a value for $m^{\mathrm{center}}$. If we are making all agents robust, we set $s^{\mathrm{center}} = s_{\mathrm{att}}^{\mathrm{center}}, s^{\mathrm{leaf }j} = s_{\mathrm{att}}^{\mathrm{leaf}}$, and set a value $m = m^{\mathrm{center}} = m^{\mathrm{leaf }j}$.

\paragraph*{Central node resilience under leaf perturbations}
We first make only the center robust.
To assess resilience, we consider the adversarial policy
\begin{equation}\label{eq:network:adversary}
    \pi^{\mathrm{leaf} j,\mathrm{adv}}_t(x_t)
    =
    \begin{cases}
        d^{\mathrm{center},\mathrm{leaf} j}_{t}(x_t) & j\in\{2,3\},
        \\
        \pi_t^{\mathrm{Nom},\mathrm{leaf} j}(x_t) & j\notin\{2,3\},
    \end{cases}
\end{equation}
where $\pi_t^{\mathrm{Nom},\mathrm{leaf} j}(\cdot)$ is the equilibrium policy of the leaf when $m^{\mathrm{center}} \to \infty$ and $d^{\mathrm{center},\mathrm{leaf} j}_{t}$ is the policy of the adversary of the center node when $m^{\mathrm{center}} = 0.05$.
We evaluate four scenarios: Nominal ($m^{\mathrm{center}} \to \infty$) or strategically robust control at the center ($m^{\mathrm{center}} = 0.05$), each with nominal or adversarial leaves~2 and~3. Table~\ref{tab:results} reports the center's terminal state and realized physical cost. Under perturbation, the robust policy stays closer to the target at lower cost.

\begin{table}[!t]
\centering
\begin{minipage}{0.30\columnwidth}
\centering
\vspace{-0.25\baselineskip}
\begin{tikzpicture}[
    scale=0.72,
    every node/.style={circle, draw, fill=gray!15, minimum size=6.5pt, inner sep=1.8pt},
    every path/.style={line width=1.1pt}
]
    \node (c) at (0,0) {};
    \node (n2) at (90:1.35) {};
    \node (n3) at (0:1.35) {};
    \node (n4) at (-90:1.35) {};
    \node (n5) at (180:1.35) {};
    \draw (c) -- (n2);
    \draw (c) -- (n3);
    \draw (c) -- (n4);
    \draw (c) -- (n5);
\end{tikzpicture}

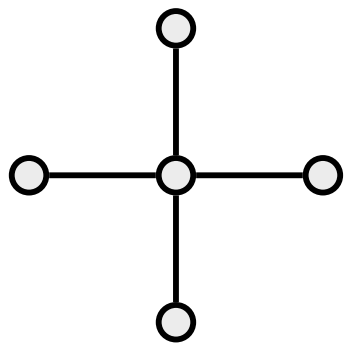
\captionof{figure}{Star graph with $N=5$ nodes.}\label{fig:network:star}
\end{minipage}
\hfill
\begin{minipage}{0.66\columnwidth}
\centering
\caption{Central-node performance under adversarial perturbation}
\vspace{-1.5ex}
\label{tab:results}
\footnotesize%
\begin{tabular}{lcc}
\toprule
Scenario & \shortstack{Terminal state\\$x_T^1$} & \shortstack{Central node\\ realized cost} \\
\midrule
Nom. & 4.09 & 2.45 \\
SR & 4.76 & 4.51 \\
Nom. (adv.) & 2.22 & 19.70 \\
SR (adv.) & 4.21 & 9.18 \\
\bottomrule
\end{tabular}
\end{minipage}
\vspace{-10pt}
\end{table}

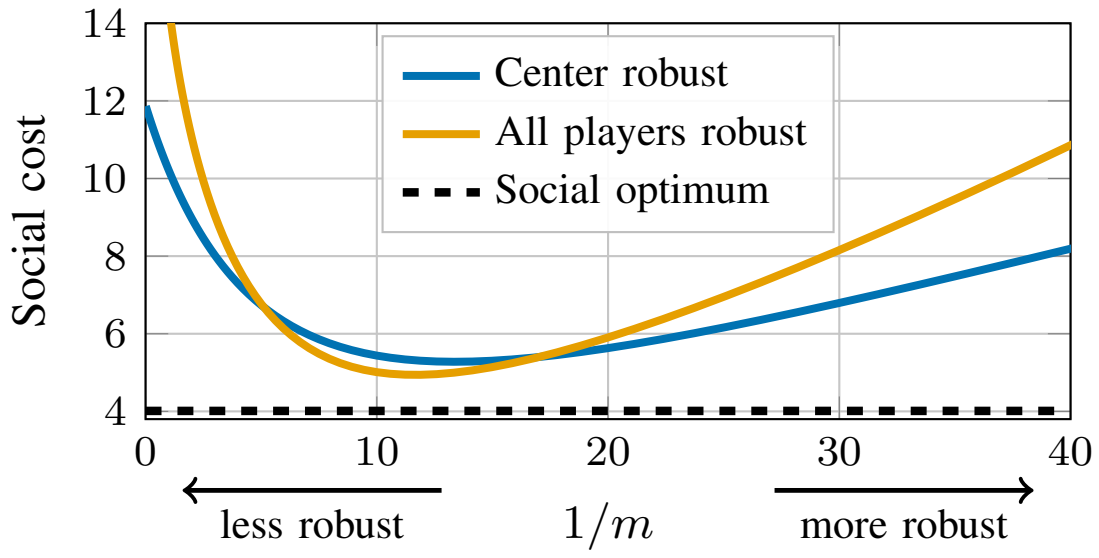
\begin{figure}[!t]
    \centering
    \scalebox{\figscale}{\input{figures/free_lunch_plot_network}}%
    \par\vspace{-2.0ex}
    \caption{Social cost, defined as the sum of the agents' realized physical costs, as a function of the strategic-robustness level $1/m$ when only the center node is robust (blue) and when everyone is (yellow).}
    \label{fig:network:free-lunch}
    \vspace{-1.75ex}
\end{figure}

\paragraph*{Robustness induces coordination}
We now repeat the social-cost study above.
We compare the case in which only the center is robust with the case of all players being robust.
As shown in Figure~\ref{fig:network:free-lunch}, the social cost again decreases over part of the robustness range, providing a second numerical instance of the free-lunch effect.

\subsection{Open-loop experiments}
We first consider a vehicle platoon, which trades a target displacement against a one-sided penalty for closing the gap to the vehicle ahead, and then a capacity accumulation game, in which competing firms invest under budget and per-stage limits toward a contest or Cournot reward at the horizon. In both, at the parameters reported below, the action sets are convex and the augmented objective is convex in $u^i$ and concave in $d^i$. The experiments are deliberately set outside of the realm of Theorem~\ref{thm:monotonicity}; we do not have analytical guarantees for monotonicity on the whole space. In practice, we still observe fast convergence. Every instance of the plain projected extragradient iteration~\cite{korpelevich1976extragradient} converges to a point with vanishing projected first-order residual. On convex--concave problems over convex sets, these conditions are globally sufficient, so the point satisfies~\eqref{eq:aug:ol} and Lemma~\ref{lemma:equivalence:ol}(a) certifies its player component as a strategically robust dynamic equilibrium.

\begin{figure}[!t]
    \centering
    \vspace{-4.5pt}
    \scalebox{0.88}{\input{figures/platoon_illustration}}
    \vspace{-14.0pt}
    \caption{Platoon of $N=4$ vehicles in one lane. Vehicle~$i$ chooses its displacement $u_t^i$; follower~$i$ keeps the gap $g_t^i$ to the vehicle ahead, with $g_{t+1}^i=g_t^i-u_t^i+u_t^{i-1}$.}
    \label{fig:platoon:illustration}
    \par\vspace{0pt}
    \scalebox{\figscale}{\input{figures/platoon_trajectories}}
    \par\vspace{-2.5ex}%
    \caption{Unperturbed equilibrium gaps $g^i_t$ between each follower and the vehicle ahead for two
    robustness levels (solid: $m=160$; dashed: $m=1000$). (a) Tight initial spacing
    ($g^i_0 = 5$): the safety penalty is active (shaded), and the robust followers open a
    margin above the desired gap $g_{\mathrm{desired}}=20$, while the near-Nash followers
    approach it from below. (b) Loose spacing ($g^i_0 = 20$): the penalty never binds, and
    the two robustness levels coincide.}
    \label{fig:platoon:trajectories}
    \vspace{-2.5ex}
\end{figure}

\begin{figure*}[!t]
    \centering
    \begin{minipage}[t]{0.37\textwidth}
        \vspace{0pt}
        \centering
        \scalebox{0.79}{\input{figures/tullock_plot}}
        \par\vspace{-10pt}%
        \caption{Equilibrium trajectories for the Tullock contest. Firms conserve their budgets early and saturate production constraints late; heightened robustness (smaller $m$) suppresses investment, as firms anticipate the adversary flooding the market.}
        \label{fig:tullock_controls}
    \end{minipage}%
    \hfill
    \begin{minipage}[t]{0.37\textwidth}
        \vspace{0pt}
        \centering
        \scalebox{0.79}{\input{figures/cournot_plot}}
        \par\vspace{-10pt}%
        \caption{Equilibrium trajectories for the Cournot competition. Similarly to the Tullock contest, players saturate production constraints late; heightened robustness (smaller $m$) suppresses investment.}
        \label{fig:cournot_controls}
    \end{minipage}%
    \hfill
    \begin{minipage}[t]{0.21\textwidth}
        \vspace{0pt}
        \centering
        \setlength{\tabcolsep}{4pt}%
        \captionof{table}{Realized player return versus robustness level $m$}
        \label{tab:roi:both}
        \vspace{-1.5ex}%
        \footnotesize%
        \begin{tabular}{cc@{\hskip 1.2em}cc}
            \toprule
            \multicolumn{2}{c}{\shortstack{Tullock\\contest}} & \multicolumn{2}{c}{\shortstack{Cournot\\competition}}\\
            \cmidrule(r){1-2}\cmidrule(l){3-4}
            $m$ & Return & $m$ & Return\\
            \midrule
            0.4 & 3.692 & 5.0 & 46.31\\
            0.2 & 4.305 & 3.0 & 50.55\\
            0.1 & 4.840 & 2.0 & 52.97\\
            0.05 & 5.456 & 1.0 & 48.48\\
            \bottomrule
        \end{tabular}
    \end{minipage}
    \vspace{-15pt}%
\end{figure*}

\subsubsection{Platoon game}
Game-theoretic formulations of vehicle platooning have been studied, e.g., in~\cite{gattami2011establishing}. We consider $N$ vehicles that travel in one lane over a horizon $T$; vehicle~1 leads and vehicles $2,\ldots,N$ follow in order. Vehicle $i$ chooses its displacement $u_t^i$. If $g_t^i$ is the gap between follower $i$ and the preceding vehicle, then $g_{t+1}^i=g_{t}^i-u_t^i+u_t^{i-1}$ for $i=2,\ldots,N$. Stacking the gaps recovers the linear dynamics~\eqref{eq:system_dynamics}.
Each follower balances a target displacement and safe spacing through
$l_t^i(g_t^i, u_t^i)=r(u_t^i-u_{\mathrm{desired}})^2+q\max\{0,g_{\mathrm{desired}}-g_t^i\}^2$; the leader has only the first term. We set $p_t^i(d_t^i, u_t^{-i})=m\|d_t^i-u_t^{-i}\|^2$ and take $N=4$, $T=20$, $r=20$, $q=1$, $g_{\mathrm{desired}}=20$, and $u_{\mathrm{desired}}=5$. For $m$ above the threshold $m^\star\approx154$, the augmented game is well posed---strictly convex in each $u^i$ and strictly concave in each $d^i$---and projected extragradient converged in all tested instances. A schematic representation is given in Figure~\ref{fig:platoon:illustration}. Figure~\ref{fig:platoon:trajectories} shows the unperturbed equilibrium gaps for two robustness levels and two initial spacings. When the platoon starts crammed (initial gap $5$), the safety penalty is active, and the robustness level visibly reshapes how the followers approach the desired spacing: the robust followers ($m=160$) settle above it, with the margin growing toward the back of the platoon, whereas the near-Nash followers ($m=1000$) approach it from below. When it starts well spaced (initial gap $20$), the penalty never binds, and the gaps are independent of the robustness level. With the same robustness parameters, robustness thus acts only when the constraint is active.

\subsubsection{Capacity accumulation game}
Closely related dynamic capacity-investment and accumulated-stock Tullock games have been studied, e.g., in~\cite{grossmann2011transitional}. In this capacity accumulation game, $N$ firms acquire inventory over a horizon $T$, with $u_t^i \in \reals$ the quantity acquired by firm $i$ at time $t$. The accumulated capacity of each firm evolves according to the decoupled linear dynamics $x_{t+1}^i = (1 -\delta) x_t^i + u_t^i$, where $\delta \in (0, 1]$ represents a depreciation factor. Stacking the individual states embeds these dynamics into~\eqref{eq:system_dynamics}.
Each firm faces a per-stage acquisition limit and a horizon budget, giving the convex action set
$
\mathcal{U}^i = \big\{ u^i \in \reals^T \,\big|\, u_t^i \in [0, \bar{u}_t^i], \,\sum_{t=0}^{T-1} u_t^i \le B^i \big\}.
$
The running cost is $l_t^i(x_t,u_t^i)=c^iu_t^i+\varepsilon(u_t^i)^2$, where $c^i>0$ and $\varepsilon>0$ models increasing marginal acquisition cost and ensures strict convexity. We set $p_t^i(d_t^i,u_t^{-i})=m\|d_t^i-u_t^{-i}\|^2$.
At time $T$, firms receive competitive rewards based on final capacities. We consider (1) the \emph{market-share (Tullock) loss}~\cite{tullock1980efficient}, $l_T^i(x_T)=-V\,x_T^i/\sum_jx_T^j$, where $V>0$ is the value of the contested market, and (2) the \emph{Cournot loss}, $l_T^i(x_T)=-x_T^i(D-\sum_jx_T^j)$, where $D>0$ is the demand intercept. We have that $J^i_{\mathrm{aug}}$ is strongly convex in $u^i$ and strongly concave in $d^i$; after we eliminate the state, we get a strategically robust dynamic equilibrium. The open-loop model is natural here because firms typically choose aggregate production plans without observing competitors' instantaneous inventories.
For the Tullock contest, we consider $N=4$ firms over a horizon $T=10$. We set initial capacities $x_0^i=1$, unit costs $c^i=1$, market value $V=30$, total budgets $B^i=8$, stage limits $\bar{u}_t^i=1$, depreciation $\delta=0.05$, and regularization $\varepsilon=10^{-6}$. We evaluate the game across varying robustness levels $m\in\{0.05, 0.1, 0.2, 0.4\}$. The resulting optimal investments and profits are illustrated in Figure~\ref{fig:tullock_controls} and Table~\ref{tab:roi:both}.
For the Cournot formulation, we maintain the same baseline setup ($N=4, T=10, x_0^i=1, \delta=0.05, \varepsilon=10^{-6}$) but adjust the economic parameters to $c^i=1$, budgets $B^i=9$, stage limits $\bar{u}_t^i=1$, and market intercept $D=30$. We evaluate $m\in\{1.0,2.0,3.0,5.0\}$. This range differs from the Tullock experiment because the two terminal losses have very different conditioning: the Tullock ratio loss is particularly poorly conditioned, so the corresponding numerical risk profiles occur on different $m$-scales. The results are shown in Figure~\ref{fig:cournot_controls}.
Table~\ref{tab:roi:both} reports the equilibrium payoffs (the negative of the realized economic loss).
Both scenarios exhibit the ``free-lunch'' phenomenon discussed in the closed-loop experiments section. In the symmetric Tullock contest, greater robustness improves individual profits over the tested range. In the Cournot game, profits are maximized at an intermediate robustness level ($m=2.0$). Curtailing production to hedge against strategic uncertainty, firms implicitly restrict aggregate supply, saving on production costs and raising net payoffs.

\section{Conclusion}\label{sec:conclusion}

We introduced strategically robust dynamic equilibria, a novel solution concept for multi-agent control that interpolates between game-theoretic control (Nash equilibrium) and robust control ($\mathcal H_\infty$).
We established open-loop existence and strong-monotonicity conditions, finite-horizon closed-loop Riccati characterizations, and a local $N$-player infinite-horizon contraction result near the decoupled $\mathcal H_\infty$ gains.
Besides providing robustness to uncertain behavior of the other agents, strategically robust dynamic equilibria are also more amenable to computation.
The numerical examples show resilience to policy misspecification and uncover a ``free-lunch'' phenomenon, in that, by acting more robustly, agents improve their nominal individual or aggregate performance.

\section{Proofs}\label{sec:proofs}

Section~\ref{sec:proof_general} isolates the equivalence between the two-stage
and the augmented formulations; Sections~\ref{sec:proofs_ol} and~\ref{sec:proofs_cl} specialize it
and prove the remaining results.

\subsection{A general equivalence principle}\label{sec:proof_general}

Lemmas~\ref{lemma:equivalence:ol} and~\ref{lemma:equivalence:cl} are two
instances of a single statement. Let $\mathcal P^i$ and $\mathcal S^i$ be nonempty action sets for player $i$ and their adversary, let
$J_{\mathrm{aug}}^i(\pi^i,\sigma^i,\pi^{-i})$ denote the augmented cost of
player $i$, and let
$J^i(\pi^i,\pi^{-i})\coloneqq\sup_{\sigma^i\in\mathcal S^i}J_{\mathrm{aug}}^i(\pi^i,\sigma^i,\pi^{-i})$
be the associated robust cost, as in~\eqref{eq:robust_cost}. We call
$(\bar\pi^1,\ldots,\bar\pi^N)$ a strategically robust dynamic equilibrium if
\begin{equation}\label{eq:sre:gen}
    J^i(\bar\pi^i,\bar\pi^{-i})
    \leq
    J^i(\pi^i,\bar\pi^{-i})
    \qquad
    \forall \pi^i\in\mathcal P^i,
\end{equation}
and $(\bar\pi^1,\bar\sigma^{1},\ldots,\bar\pi^N,\bar\sigma^{N})$ an
equilibrium of the augmented game if
\begin{equation}\label{eq:aug:gen}
    J_\mathrm{aug}^i(\bar \pi^{i}, \sigma^i, \bar\pi^{-i})
    \leq
    J_\mathrm{aug}^i(\bar \pi^{i}, \bar\sigma^i, \bar\pi^{-i})
    \leq
    J_\mathrm{aug}^i(\pi^i, \bar\sigma^i, \bar\pi^{-i})
\end{equation}
for all $\pi^i\in\mathcal P^i$ and $\sigma^i\in\mathcal S^i$; both conditions
are imposed for every player $i$. Definitions~\ref{def:sre:ol}
and~\ref{def:aug:ol} are the case
$(\mathcal P^i,\mathcal S^i)=(\mathcal U^i,\mathcal U^{-i})$, whereas
Definitions~\ref{def:sre:cl} and~\ref{def:aug:cl} are the case
$(\mathcal P^i,\mathcal S^i)=(\Pi^i,\Sigma^i)$ of state-feedback policies.

\begin{lemma}[general equivalence]\label{lemma:equivalence:gen}
    For each player $i$ and fixed opponents' policies $\bar\pi^{-i}$, consider the
    zero-sum game in which player $i$ minimizes $J^i_{\mathrm{aug}}$ over $\pi^i$
    and their adversary maximizes over $\sigma^i$. Then the following hold.
    \begin{enumerate}
        \item[\textnormal{(a)}] \emph{Necessity.}
        If $(\bar\pi^1,\bar\sigma^{1},\ldots,\bar\pi^N,\bar\sigma^{N})$ is an
        equilibrium of the augmented game, then $(\bar\pi^1,\ldots,\bar\pi^N)$ is
        a strategically robust dynamic equilibrium and each $\bar\sigma^i$ is a
        worst-case policy in~\eqref{eq:robust_cost}, i.e., a best response to
        $\bar\pi^i$.
        \item[\textnormal{(b)}] \emph{Sufficiency.} Suppose
        $(\bar\pi^1,\ldots,\bar\pi^N)$ is a strategically robust dynamic
        equilibrium and that, for each $i$, the zero-sum game
        $J^i_{\mathrm{aug}}(\cdot,\cdot,\bar\pi^{-i})$ admits a saddle point $(\pi^{i,\star},\sigma^{i,\star})$. Then, taking each $\bar\sigma^i$ to be
        the adversary's component $\sigma^{i,\star}$ of such a saddle point, the
        tuple $(\bar\pi^1,\bar\sigma^{1},\ldots,\bar\pi^N,\bar\sigma^{N})$ is an
        equilibrium of the augmented game; this $\bar\sigma^i$ is moreover one of
        the worst-case policies in~\eqref{eq:robust_cost}, i.e., a best response to
        $\bar\pi^i$.
    \end{enumerate}
\end{lemma}

\begin{proof}
    Fix $\bar\pi^{-i}$; player $i$'s problem is a zero-sum game against its
    adversary with objective $J^i_{\mathrm{aug}}$. By definition, the pair
    $(\bar\pi^i,\bar\sigma^i)$ satisfies the augmented-equilibrium
    condition~\eqref{eq:aug:gen} exactly when it is a saddle point of
    $J^i_{\mathrm{aug}}(\cdot,\cdot,\bar\pi^{-i})$.

    \emph{(a)} If $(\bar\pi^i,\bar\sigma^i)$
    satisfies~\eqref{eq:aug:gen}, then $\bar\sigma^i$ attains
    $\sup_{\sigma^i}J^i_{\mathrm{aug}}(\bar\pi^i,\sigma^i,\bar\pi^{-i})$. Hence, for
    any alternative $\pi^i$,
    \begin{equation*}
    \begin{aligned}
        \sup_{\sigma^i} J_{\mathrm{aug}}^i(\bar\pi^i,\sigma^i,\bar\pi^{-i})
        &= J_{\mathrm{aug}}^i(\bar\pi^i,\bar\sigma^i,\bar\pi^{-i})
        \le J_{\mathrm{aug}}^i(\pi^i,\bar\sigma^i,\bar\pi^{-i}) \\
        &\le \sup_{\sigma^i} J_{\mathrm{aug}}^i(\pi^i,\sigma^i,\bar\pi^{-i}),
    \end{aligned}
    \end{equation*}
    which is~\eqref{eq:sre:gen}.

    \emph{(b)} Let $(\bar\pi^1,\dots,\bar\pi^N)$ be a strategically robust dynamic
    equilibrium and fix $i$. Let $(\pi^{i,\star},\sigma^{i,\star})$ be a saddle
    point of $J^i_{\mathrm{aug}}(\cdot,\cdot,\bar\pi^{-i})$, with value
    $V^i\coloneqq J^i_{\mathrm{aug}}(\pi^{i,\star},\sigma^{i,\star},\bar\pi^{-i})$;
    set $\bar\sigma^i\coloneqq\sigma^{i,\star}$. The saddle inequalities read, for
    all $\pi^i,\sigma^i$,
    \begin{equation}\label{eq:saddle}
        J^i_{\mathrm{aug}}(\pi^{i,\star},\sigma^i,\bar\pi^{-i})
        \;\le\; V^i \;\le\;
        J^i_{\mathrm{aug}}(\pi^i,\sigma^{i,\star},\bar\pi^{-i}).
    \end{equation}
    We verify~\eqref{eq:aug:gen} for $(\bar\pi^i,\bar\sigma^i)$. The right
    inequality of~\eqref{eq:saddle}, with $\bar\sigma^i=\sigma^{i,\star}$, gives
    $J^i_{\mathrm{aug}}(\pi^i,\bar\sigma^i,\bar\pi^{-i})\ge V^i$ for every $\pi^i$;
    at $\pi^i=\bar\pi^i$ this reads
    $J^i_{\mathrm{aug}}(\bar\pi^i,\bar\sigma^i,\bar\pi^{-i})\ge V^i$. Conversely,
    the strategically robust equilibrium property~\eqref{eq:sre:gen} and the fact
    that $\sigma^{i,\star}$ maximizes against $\pi^{i,\star}$ give
    $$
        \sup_{\sigma^i} J_{\mathrm{aug}}^i(\bar\pi^i,\sigma^i,\bar\pi^{-i})
        \;\le\;
        \sup_{\sigma^i} J_{\mathrm{aug}}^i(\pi^{i,\star},\sigma^i,\bar\pi^{-i})
        \;=\; V^i,
    $$
    so in particular
    $J^i_{\mathrm{aug}}(\bar\pi^i,\bar\sigma^i,\bar\pi^{-i})\le V^i$. Combining the two,
    $$
        J^i_{\mathrm{aug}}(\bar\pi^i,\bar\sigma^i,\bar\pi^{-i}) = V^i
        = \sup_{\sigma^i} J_{\mathrm{aug}}^i(\bar\pi^i,\sigma^i,\bar\pi^{-i}),
    $$
    which is the left inequality of~\eqref{eq:aug:gen};
    $\bar\sigma^i$ attains the supremum in~\eqref{eq:robust_cost}. The right
    inequality
    of~\eqref{eq:aug:gen},
    $J^i_{\mathrm{aug}}(\bar\pi^i,\bar\sigma^i,\bar\pi^{-i})=V^i\le
    J^i_{\mathrm{aug}}(\pi^i,\bar\sigma^i,\bar\pi^{-i})$ for all $\pi^i$, is the
    right inequality of~\eqref{eq:saddle}. Hence $(\bar\pi^i,\bar\sigma^i)$ is an
    augmented-game equilibrium for player $i$; as $i$ was arbitrary, the claim
    follows.
\end{proof}

\subsection{Proofs for the open-loop setting}\label{sec:proofs_ol}

\begin{lemma}[transfer of Assumption~\ref{ass:ol}]\label{lemma:transfer}
    Let Assumption~\ref{ass:ol} hold. Then, for every player $i$ and at every
    joint profile,
    \begin{equation}\label{eq:induced}
        \nabla^2_{u^iu^i}J^i_{\mathrm{aug}}\succeq\mu_U^iI,
        \qquad
        -\nabla^2_{d^id^i}J^i_{\mathrm{aug}}\succeq\mu_D^iI,
    \end{equation}
    with the constants of~\eqref{eq:constants}. If in addition $\mu_U^i>0$ and
    $\mu_D^i>0$, then $u^i\mapsto J^i_{\mathrm{aug}}$ is strongly convex and
    coercive, and $d^i\mapsto J^i_{\mathrm{aug}}$ is strongly concave and tends
    to $-\infty$ as $\norm{d^i}\to\infty$.
\end{lemma}

\begin{proof}
Each stage cost depends only on the time-$t$ blocks of its arguments, so the
Hessians of $l^i$ and $p^i$ act blockwise and inherit the stagewise bounds of
Assumption~\ref{ass:ol}: for arbitrary $x$ and $v$,
$x^\top\nabla^2_{xx}l^ix=\sum_tx_t^\top\nabla^2_{xx}l_t^ix_t
\ge\underline q^i\sum_t\norm{x_t}^2=\underline q^i\norm{x}^2$, and likewise for
the other definite bounds, while for the cross terms
$x^\top\nabla^2_{xu^i}l^iv=\sum_tx_t^\top\nabla^2_{xu^i}l_t^iv_t
\le\kappa_u^i\sum_t\norm{x_t}\,\norm{v_t}\le\kappa_u^i\norm{x}\,\norm{v}$
by Cauchy--Schwarz.
Substituting $x=G^iu^i+G^{-i}d^i+Cx_0$ into $l^i$ and applying the chain
rule,
\begin{equation*}
\begin{aligned}
    \nabla^2_{u^iu^i}J^i_{\mathrm{aug}}
    &=(G^i)^\top\nabla^2_{xx}l^i\,G^i+\nabla^2_{u^iu^i}l^i\\
    &\quad+(G^i)^\top\nabla^2_{xu^i}l^i+(\nabla^2_{xu^i}l^i)^\top G^i;
\end{aligned}
\end{equation*}
the identity for $d^i$ is the same with $G^{-i}$ in place of $G^i$ and the
extra term $-\nabla^2_{d^id^i}p^i$. Bounding the groups of terms through the
stacked bounds just established and using $(G^i)^\top\nabla^2_{xu^i}l^i+(\nabla^2_{xu^i}l^i)^\top G^i \succeq-2\norm{(G^i)^\top\nabla^2_{xu^i}l^i} I$
gives~\eqref{eq:induced}.
When $\mu_U^i>0$, the first bound in~\eqref{eq:induced} makes
$u^i\mapsto J^i_{\mathrm{aug}}$ $\mu_U^i$-strongly convex. From any fixed
$\bar u^i$ it is bounded below by
$J^i_{\mathrm{aug}}(\bar u^i,\cdot,\cdot)
-\norm{\nabla_{u^i}J^i_{\mathrm{aug}}(\bar u^i,\cdot,\cdot)}\,r
+\tfrac12\mu_U^ir^2$ at $r=\norm{u^i-\bar u^i}$, which tends to
$+\infty$; analogously $\mu_D^i>0$ sends
$d^i\mapsto J^i_{\mathrm{aug}}$ to $-\infty$.
\end{proof}

\noindent The proofs below use~\eqref{eq:induced}, the bounds $c^{ij}$,
and the penalty properties stated in Section~\ref{sec:setting}.

\begin{proof}[Proof of Lemma~\ref{lemma:equivalence:ol}]
Apply Lemma~\ref{lemma:equivalence:gen} with $\mathcal P^i=\mathcal U^i$,
$\mathcal S^i=\mathcal U^{-i}$, and the augmented cost $J_{\mathrm{aug}}^i$ of
Section~\ref{sec:setting}.
Part~(a) is exactly part~(a) of Lemma~\ref{lemma:equivalence:gen} and requires
no assumption. For part~(b) it remains to exhibit, for each $i$ and each fixed
$\bar u^{-i}$, a saddle point of
$f\coloneqq J^i_{\mathrm{aug}}(\cdot,\cdot,\bar u^{-i})$ on
$\mathcal U^i\times\mathcal U^{-i}$.
Fix any $\hat u\in\mathcal U^i$ and $\hat d\in\mathcal U^{-i}$. The map
$f(\cdot,\hat d)$ is continuous, strictly convex, and coercive on the nonempty
closed convex set $\mathcal U^i$, so
$B\coloneqq\min_{u\in\mathcal U^i}f(u,\hat d)$ is finite and attained;
symmetrically, $A\coloneqq\max_{d\in\mathcal U^{-i}}f(\hat u,d)$ is finite and
attained. For $r>\max\{\norm{\hat u},\norm{\hat d}\}$ set
$\mathcal U_r\coloneqq\{u\in\mathcal U^i:\norm{u}\le r\}$ and
$\mathcal D_r\coloneqq\{d\in\mathcal U^{-i}:\norm{d}\le r\}$, which are nonempty,
compact, and convex. On $\mathcal U_r\times\mathcal D_r$ the function $f$ is
continuous and convex--concave, so by Sion's minimax theorem $f$ has a saddle
point $(u_r,d_r)$ there, attained by compactness.

These truncated saddle points are bounded uniformly in $r$. Evaluating the
saddle inequalities at $d=\hat d$ and at $u=\hat u$ gives
$
    B\le f(u_r,\hat d)\le f(u_r,d_r)\le f(\hat u,d_r)\le A,
$
so $u_r$ lies in the sublevel set $\{u\in\mathcal U^i:f(u,\hat d)\le A\}$ and
$d_r$ in the superlevel set $\{d\in\mathcal U^{-i}:f(\hat u,d)\ge B\}$, both
bounded by coercivity. Let $R$ bound them and fix any $r>R$. Then
$\norm{u_r}<r$, so for every $u\in\mathcal U^i$ the point
$u_\lambda\coloneqq u_r+\lambda(u-u_r)$ lies in $\mathcal U_r$ for all
sufficiently small $\lambda>0$, and convexity of $f(\cdot,d_r)$ gives
$f(u_r,d_r)\le f(u_\lambda,d_r)\le(1-\lambda)f(u_r,d_r)+\lambda f(u,d_r)$,
whence $f(u_r,d_r)\le f(u,d_r)$. The argument in $d$ is symmetric. Hence
$(u_r,d_r)$ is a saddle point of $f$ on all of
$\mathcal U^i\times\mathcal U^{-i}$. Lemma~\ref{lemma:equivalence:gen}(b)
therefore applies.
\end{proof}

\begin{proof}[Proof of Proposition~\ref{prop:ole}]
\emph{(a)} Each $\mathcal U^i$ is compact and convex by assumption, and each
$\mathcal U^{-i}$ is a finite product of such sets. Best responses exist by continuity and
compactness and are unique by the strict convexity--concavity that
$\mu_U^i,\mu_D^i>0$ gives through~\eqref{eq:induced}, so Berge's Maximum
Theorem makes them continuous.
The joint best-response map thus sends
$\prod_i(\mathcal U^i\times\mathcal U^{-i})$ continuously into itself and, by
Brouwer's Fixed-Point Theorem, has a fixed point. This is an equilibrium of the
augmented game; Lemma~\ref{lemma:equivalence:ol}(a) gives the claim.

\emph{(b)} Fix $i$. The problem is unconstrained and $J^i_{\mathrm{aug}}$ is
strictly convex in $u^i$, so $\nabla_{u^i}J^i_{\mathrm{aug}}=0$ is necessary
and sufficient for ${u^i}^\star$ to minimize
$u^i\mapsto J^i_{\mathrm{aug}}(u^i,{d^i}^\star,{u^{-i}}^\star)$; by strict concavity,
$\nabla_{d^i}J^i_{\mathrm{aug}}=0$ is necessary and sufficient for ${d^i}^\star$
to maximize $d^i\mapsto J^i_{\mathrm{aug}}({u^i}^\star,d^i,{u^{-i}}^\star)$. These are exactly~\eqref{eq:ol:foc}. Holding for every $i$, they give~\eqref{eq:aug:ol};
Lemma~\ref{lemma:equivalence:ol}(a) concludes.
\end{proof}

\begin{proof}[Proof of Theorem~\ref{thm:monotonicity}]
Order $z=(u^1,\ldots,u^N,d^1,\ldots,d^N)$ and let $H_\lambda\coloneqq\tfrac12(\nabla F_\lambda+\nabla F_\lambda^\top)$. Since $J^i_{\mathrm{aug}}$ depends on $u^j$, $j\neq i$, only through $p^i$, which does not involve $u^i$, and does not depend on $d^j$, $j\neq i$, the $u$-$u$ and $d$-$d$ parts of $H_\lambda$ are block-diagonal, with blocks $\lambda_k\nabla^2_{u^ku^k}J^k_{\mathrm{aug}}\succeq\lambda_k\mu_U^kI$ and $-\lambda_k\nabla^2_{d^kd^k}J^k_{\mathrm{aug}}\succeq\lambda_k\mu_D^kI$ by~\eqref{eq:induced}. The $(u^i,d^i)$ blocks are exact transposes with opposite signs and equal weights, so they cancel. The block pairing $u^j$ with $d^i$ is $\tfrac12\lambda_i(\nabla^2_{d^iu^j}p^i)^\top$, of norm at most $\tfrac12\lambda_ic^{ij}$. Thus, for $v=\operatorname{col}(a^1,\ldots,a^N,b^1,\ldots,b^N)$,
\begin{equation*}
    v^\top H_\lambda v\ge\sum_k\lambda_k\big(\mu_U^k\norm{a^k}^2+\mu_D^k\norm{b^k}^2\big)-\sum_{i\neq j}\lambda_ic^{ij}\norm{b^i}\norm{a^j}.
\end{equation*}
Since $p_t^i\ge0$ and $p_t^i(y,y)=0$, we have $\nabla_{d^i}p_t^i(y,y)=0$ for all $y$, and differentiating in $y$ gives $\nabla^2_{d^iu^{-i}}p_t^i=-\nabla^2_{d^id^i}p_t^i\preceq-\mu_p^iI$ at $d^i=u^{-i}$. The block of this matrix for $u^j$ contains a diagonal block of $-\nabla^2_{d^id^i}p_t^i$, so $c^{ij}\ge\mu_p^i$ for all $j\neq i$; in particular, $c_{\mathrm{in}}^k,c_{\mathrm{out}}^k>0$ and $\chi^k<\infty$ for every $k$. Let $\theta\coloneqq(\chi^{(1)}\chi^{(2)})^{1/4}>1$. We check that the weights~\eqref{eq:weights} satisfy
\begin{equation}\label{eq:edge_weights}
    4\mu_D^i\mu_U^j\lambda_j\ge\theta^2c_{\mathrm{out}}^ic_{\mathrm{in}}^j\lambda_i\qquad\text{for all }i\neq j.
\end{equation}
Write $\lambda_k=w_kc_{\mathrm{in}}^k/\mu_U^k$, with $w_k=1$ for $k\neq k^\star$ and $w_{k^\star}=\sqrt{\chi^{(1)}/\chi^{(2)}}$. Since $c_{\mathrm{in}}^j>0$, \eqref{eq:edge_weights} is equivalent to $w_j\chi^i\ge\theta^2w_i$. For $i=k^\star$ it holds with equality. For $j=k^\star$ it reads $\chi^i\ge \chi^{(2)}$, and for $i,j\neq k^\star$ it reads $\chi^i\ge\sqrt{\chi^{(1)}\chi^{(2)}}$; both hold because $\chi^i\ge \chi^{(2)}\ge\sqrt{\chi^{(1)}\chi^{(2)}}$ for $i\neq k^\star$.

For all $i\neq j$, the AM--GM inequality and~\eqref{eq:edge_weights} give
\begin{equation*}
    \lambda_ic^{ij}\norm{b^i}\norm{a^j}
    \le\frac{c^{ij}}{\theta}\Big(\frac{\lambda_j\mu_U^j}{c_{\mathrm{in}}^j}\norm{a^j}^2+\frac{\lambda_i\mu_D^i}{c_{\mathrm{out}}^i}\norm{b^i}^2\Big).
\end{equation*}
Summing over $i\neq j$ and grouping the terms by $j$ and by $i$, respectively, gives
\begin{align*}
    \sum_{i\neq j}\lambda_ic^{ij}\norm{b^i}\norm{a^j}
    \le \frac1\theta\sum_k\lambda_k\big(\mu_U^k\norm{a^k}^2+\mu_D^k\norm{b^k}^2\big),
\end{align*}
since the inner sums equal $c_{\mathrm{in}}^j$ and $c_{\mathrm{out}}^i$ by~\eqref{eq:cycle_constants}. Combined with the first estimate, this yields
\begin{equation*}
    v^\top H_\lambda v\ge\Big(1-\frac1\theta\Big)\sum_k\lambda_k\big(\mu_U^k\norm{a^k}^2+\mu_D^k\norm{b^k}^2\big)
\end{equation*}
uniformly in $z$. Since $\theta>1$, $F_\lambda$ is strongly monotone.

If $c_{\mathrm{out}}^ic_{\mathrm{in}}^j<4\mu_D^i\mu_U^j$ for all $i\neq j$, \eqref{eq:edge_weights} also holds for $\lambda\equiv1$ and some $\theta>1$, and the same argument applies.
\end{proof}

\subsection{Proofs for the closed-loop setting}\label{sec:proofs_cl}

\begin{proof}[Proof of Lemma~\ref{lemma:equivalence:cl}]
Apply Lemma~\ref{lemma:equivalence:gen} with $\mathcal P^i=\Pi^i$ and
$\mathcal S^i=\Sigma^i$, the closed-loop policy sets of
Section~\ref{sec:closed_loop}. Part~(a) is part~(a) of the general
lemma, while part~(b) follows from its sufficiency direction under the
saddle-point hypothesis.
\end{proof}

\begin{proof}[Proof of Theorem~\ref{thm:sre}]
We construct a candidate equilibrium of the augmented game by backward induction, verify the saddle inequalities~\eqref{eq:aug:cl}, and prove subgame-perfect uniqueness.

\paragraph*{Step 1}
We compute a candidate equilibrium using backward induction
\begin{multline}\label{eq:proof:backward_dp}
    \!\!\!\!
    V_t^i(x) = \min_{u_t^i}\max_{d_t^i}
    \norm{x}^2_{Q_t^i} + \norm{u_t^i}^2_{R_t^i} - \norm{d_t^i}^2_{S_t^i}
    -\norm{d_t^i-u_t^{-i}}_{M_t^i}^2
    \\
    +
    V_{t+1}^i\left(A_tx + B_t^iu_t^i + B_t^{-i} d_t^i\right).
\end{multline}
We claim that $V_t^i(x)=x^\top P_t^ix$, with $P_t^i$ given by
\eqref{eq:sre_riccati}. The claim follows by backward induction. The terminal
case is immediate. Suppose $V_{t+1}^i(x)=x^\top P_{t+1}^ix$. Then the
stage-$t$ saddle problem~\eqref{eq:proof:backward_dp} gives
\begin{align*}
    0
    &=
    R_t^iu_t^i + (B_t^i)^\top P_{t+1}^i\left(A_tx_t + B_t^iu_t^i + B_t^{-i} d_t^i\right)
    \\
    0
    &=
    -S^i_td_t^i - M_t^i(d^i_t-u_t^{-i})
    \\
    &\qquad\qquad
    + (B_t^{-i})^\top P^i_{t+1}\left(A_tx_t + B_t^iu_t^i + B_t^{-i} d_t^i\right).
\end{align*}
The assumptions $H_t^{u^i,u^i}\succ0$ and
$H_t^{d^i,d^i}\succ0$ make this objective strongly convex in $u_t^i$ and
strongly concave in $d_t^i$, so the lower and upper values coincide and the
saddle point is unique for fixed $u_t^{-i}$. Stacking the first-order
conditions for all players and adversaries gives the linear system
\eqref{eq:sre_policies}; invertibility of $H_t$ therefore yields gains
$u_t^i=K_t^ix_t$ and $d_t^i=L_t^ix_t$. Substitution into the value function
gives~\eqref{eq:sre_riccati} and completes the induction.
We call this candidate equilibrium $(\bar\pi^1,\bar\sigma^1,\ldots,\bar\pi^N,\bar\sigma^N)$.

\paragraph*{Step 2}
We show that the candidate of Step~1 satisfies~\eqref{eq:aug:cl}.

For the right inequality (optimality of $\bar\pi^i$), fix $d_t^i = L_t^i x_t$ and $u_t^j = K_t^j x_t$ for all $j\neq i$.
Player~$i$ then faces a standard LQR with dynamics
$
    x_{t+1} = \tilde A_t^i x_t + B_t^i u_t^i,
$
where
$
    \tilde A_t^i \coloneqq A_t + B_t^{-i} L_t^i,
$
and cost matrices $\tilde R_t^i = R_t^i$, $\tilde Q_T^i = Q_T^i$, and
$
    \tilde Q_t^i
    = Q_t^i
    - (L_t^i)^\top S_t^i L_t^i
    - (L_t^i - K_t^{-i})^\top M_t^i (L_t^i - K_t^{-i}).
$
The LQR gain for this reduced problem is
$
    \tilde K_t^i
    =
    -\bigl(R_t^i + (B_t^i)^\top \tilde P_{t+1}^i B_t^i\bigr)^{-1} (B_t^i)^\top \tilde P_{t+1}^i \tilde A_t^i.
$
We claim $\tilde K_t^i = K_t^i$ and $\tilde P_t^i = P_t^i$ for all $t$, which we verify by backward induction.
The base case $\tilde P_T^i = Q_T^i = P_T^i$ is immediate.
Suppose $\tilde P_{t+1}^i = P_{t+1}^i$.
The $u^i$-row of the first-order conditions in Step~1 reads
$H_t^{u^i,u^i} K_t^i+H_t^{u^i,d^i}L_t^i =G_t^{u^i}$, which, by direct inspection,
gives
$
    K_t^i = \tilde K_t^i.
$
Since the closed-loop matrices coincide, $\tilde A_t^i + B_t^i \tilde K_t^i = A_t + B_t^i K_t^i + B_t^{-i} L_t^i$, the LQR Riccati recursion propagates $\tilde P_t^i = P_t^i$.
Moreover, the reduced LQR has a unique finite minimizer at every stage because $R_t^i + (B_t^i)^\top P_{t+1}^i B_t^i = H_t^{u^i,u^i} \succ 0$ by assumption.
Thus, $\bar\pi^i$ is the unique optimal policy, proving the right
inequality in~\eqref{eq:aug:cl}. Fixing instead
$u_t^i=K_t^ix_t$ and $u_t^{-i}=K_t^{-i}x_t$, the adversary faces the
corresponding strictly concave LQ maximization problem. The $d^i$-row of
\eqref{eq:sre_policies} and the same backward-induction argument identify
$L_t^i$ and $P_t^i$ as its unique optimizer and value matrix, proving the left
inequality. Hence the candidate is an equilibrium of the augmented game, and
Lemma~\ref{lemma:equivalence:cl}(a) shows that its player component is a
closed-loop strategically robust dynamic equilibrium.

\paragraph*{Step 3}
For uniqueness, suppose there is another subgame-perfect strategically robust dynamic equilibrium, and let $t'$ be
the last time at which one of its player policies differs from the candidate
above at some admissible history. Since all player policies agree at later
times, the robust continuation values at $t'+1$ coincide and are quadratic
functions of the current state; hence, the matrices $H_{t'}$ and $G_{t'}^u, G_{t'}^{d}$ are
the same as in Step~1. At every admissible history, the stage-$t'$ minimax
problem is strongly convex--concave, so each player's equilibrium action and
unique worst-case response satisfy the saddle-point first-order conditions.
Stacking these conditions gives~\eqref{eq:sre_policies}, whose solution is
unique by invertibility of $H_{t'}$ and depends only on $x_{t'}$. Thus, all
player policies must also coincide at $t'$, a contradiction.
\end{proof}

\begin{proof}[Proof of Proposition~\ref{prop:invertibility}]
Let $D_\lambda$ be the block-diagonal matrix that scales both the
$u^i$- and $d^i$-rows by the weight $\lambda_i$ of~\eqref{eq:weights},
evaluated with the constants of the proposition.
The symmetric part of $D_\lambda H_t$ has diagonal blocks
$\lambda_iH_t^{u^i,u^i}$ and $\lambda_iH_t^{d^i,d^i}$.
The same-player off-diagonal blocks cancel, while the block pairing
$u^j$ with $d^i$, $j\neq i$, is
$-\tfrac12\lambda_i(M_t^{i,j})^\top$.
These are the block structure and bounds used in the proof of
Theorem~\ref{thm:monotonicity}, so~\eqref{eq:monotonicity_condition}
makes this symmetric part positive definite.
Hence $D_\lambda H_t$ is nonsingular, and since $D_\lambda$ is
nonsingular, so is $H_t$.
\end{proof}

\noindent\textbf{Notation for the infinite horizon.}
Fix player $i$ and opponents' gains $K^{-i}$, so that
$u_t^{-i}=K^{-i}x_t$, and center the adversary by
$
\widetilde d_t^i\coloneqq
 d_t^i-\frac{m^i}{s^i+m^i}K^{-i}x_t.
$
For fixed $K^{-i}$, player $i$'s fixed-opponent best-response problem is exactly
\begin{equation}\label{eq:inf:centered_problem}
\inf_{\pi^i}\sup_{\widetilde\sigma^i}
\sum_{t=0}^{\infty}
\norm{x_t}_{\bar Q^i(K^{-i})}^2+\norm{u_t^i}_{R^i}^2
-(s^i+m^i)\norm{\widetilde d_t^i}^2,
\end{equation}
subject to
$x_{t+1}=\bar A^i(K^{-i})x_t+B^iu_t^i+B^{-i}\widetilde d_t^i$, where
\begin{align}\label{eq:inf:centered}
\bar A^i(K^{-i})&\coloneqq
A+\tfrac{m^i}{s^i+m^i}B^{-i}K^{-i},\nonumber\\
\bar Q^i(K^{-i})&\coloneqq
Q^i-\tfrac{s^im^i}{s^i+m^i}(K^{-i})^\top K^{-i}.
\end{align}
Thus~\eqref{eq:inf:centered_problem} is a time-invariant
$\mathcal H_\infty$ state-feedback problem with state weight $\bar Q^i$,
control weight $R^i$, and disturbance weight $s^i+m^i$. Whenever it admits a
finite-value stabilizing solution, write its value as $x^\top P^ix$, where
$P^i=P^i(K^{-i})$ solves
\begin{align}\label{eq:inf:dare}
  P^i&=\bar Q^i+(\bar A^i)^\top P^i\big(I+G^iP^i\big)^{-1}\bar A^i, \nonumber \\
  G^i&\coloneqq B^i(R^i)^{-1}(B^i)^\top
  -\frac{1}{s^i+m^i}B^{-i}(B^{-i})^\top.
\end{align}
We write
$\mathcal T_{m^i}^i(\cdot;K^{-i})$ for the right-hand side
of~\eqref{eq:inf:dare}.
The gains $u^i=\mathrm{BR}^ix$ and $d^i=L^ix$ solve
{\setlength{\arraycolsep}{2pt}
\begin{align}
  &\underbrace{\begin{bmatrix}
    R^i+(B^i)^\top P^iB^i & (B^i)^\top P^iB^{-i}\\[2pt]
    -(B^{-i})^\top P^iB^i & (s^i+m^i)I-(B^{-i})^\top P^iB^{-i}
  \end{bmatrix}}_{=:\,\mathcal H^i(P^i)}
  \begin{bmatrix}\mathrm{BR}^i\\[1pt]L^i\end{bmatrix}
  =\nonumber\\
  &\begin{bmatrix}-(B^i)^\top P^iA\\[1pt](B^{-i})^\top P^iA+m^iK^{-i}\end{bmatrix}.
  \label{eq:inf:saddle}
\end{align}}
The resulting gains generate the
fictitious saddle closed loop
\begin{align}\label{eq:inf:fictitious_dynamics}
F^i(K^{-i})
&\coloneqq A+B^i\mathrm{BR}^i(K^{-i})+B^{-i}L^i(K^{-i})\nonumber\\
&=\bar A^i(K^{-i})+B^i\mathrm{BR}^i(K^{-i})+B^{-i}\widetilde L^i(K^{-i})\nonumber\\
&=\big(I+G^iP^i(K^{-i})\big)^{-1}\bar A^i(K^{-i}),
\end{align}
with $\widetilde L^i=L^i-\frac{m^i}{s^i+m^i}K^{-i}$.
At $m^i=0$, \eqref{eq:inf:centered} gives $\bar A^i=A$ and $\bar Q^i=Q^i$, so~\eqref{eq:inf:centered_problem} is the standalone $\mathcal H_\infty$ problem of Assumption~\ref{ass:feas} with disturbance weight $s^i$; standard theory gives $s^iI-(B^{-i})^\top P_\infty^iB^{-i}\succ0$.

\begin{proof}[Proof of Proposition~\ref{prop:inf_verify}]
Fix player~$i$ and the opponents' gains $K^{-i}$. By~\eqref{eq:complete_square}, the change of variables $\widetilde d_t^i=d_t^i-\frac{m^i}{s^i+m^i}K^{-i}x_t$ turns player~$i$'s robust best-response problem into~\eqref{eq:inf:centered_problem}, and by~\eqref{eq:inf:verify_domain} we have $\bar Q^i\succ0$. Condition~(i) makes the one-step objective of~\eqref{eq:inf:centered_problem} with value $x^\top P^ix$ strictly convex in $u^i$ and strictly concave in $\widetilde d^i$, since $R^i+(B^i)^\top P^iB^i$ and $(s^i+m^i)I-(B^{-i})^\top P^iB^{-i}$ are the diagonal blocks of the symmetric part of $\mathcal H^i(P^i)$. Its linear system, after substituting $L^i=\widetilde L^i+\frac{m^i}{s^i+m^i}K^{-i}$, is the first-order condition of this one-step problem, so $K^i$ and $\widetilde L^i$ are its unique saddle gains. By~\eqref{eq:complete_square}, \eqref{eq:inf:care} is the Riccati identity of~\eqref{eq:inf:centered_problem} along these gains, and by~(ii) their closed loop $A^i_{\mathrm{CL}}=\bar A^i+B^iK^i+B^{-i}\widetilde L^i$ is Schur. These are the standard verification conditions for~\eqref{eq:inf:centered_problem}~\cite[Ch.~3]{bacsar2008h}, which therefore identify $x^\top P^ix$ as its value, $K^i$ as an optimal controller, and $\widetilde L^i$ as a worst-case disturbance against $K^i$. Reverting the change of variables, $L^i$ is worst-case against $K^i$, and $K^i=\mathrm{BR}^i(K^{-i})$. Since this holds for every $i$, $K=\mathrm{BR}(K)$, and hence $K$ is a stationary linear strategically robust dynamic equilibrium.
\end{proof}

\begin{proof}[Proof of Theorem~\ref{thm:inf}]
We divide the proof into three steps.

\emph{Step 1:}
We prove part~(a). At $m^i=0$, player~$i$'s problem is the decoupled $\mathcal H_\infty$ problem of Assumption~\ref{ass:feas}, whose solution is known; we continue this solution to small $m^i>0$.

To start, we show that the implicit-function theorem applies at $m^i=0$.
Fix player $i$ and, in this proof, display the dependence on $m^i$ explicitly.
Set
$$
\Phi^i(P,K^{-i},m^i)
\coloneqq P-\mathcal T_{m^i}^i(P;K^{-i}).
$$
At $m^i=0$, \eqref{eq:inf:centered} gives $\bar A^i=A$ and $\bar Q^i=Q^i$ for every $K^{-i}$, so $\mathcal T_0^i(\cdot;K^{-i})$ is the Riccati map of the decoupled $\mathcal H_\infty$ problem of Assumption~\ref{ass:feas} and $\Phi^i(P_\infty^i,K^{-i},0)=0$.
Moreover, differentiating the Riccati map with respect to $P$ gives
$
D_P\Phi^i(P,K^{-i},m^i)[\Delta]
=\Delta-F(P,K^{-i},m^i)^\top\Delta F(P,K^{-i},m^i)
$
where
$
F(P,K^{-i},m^i)
\coloneqq\bigl(I+G^i(m^i)P\bigr)^{-1}\bar A^i(K^{-i},m^i).
$
By~\eqref{eq:inf:fictitious_dynamics} at $m^i=0$, $F(P_\infty^i,K^{-i},0)=F_\infty^i=A+B^iK_\infty^i+B^{-i}L_\infty^i$, which is Schur by Assumption~\ref{ass:feas}.
Hence the Lyapunov operator
$
\Delta\mapsto
\Delta-(F_\infty^i)^\top\Delta F_\infty^i
$
is invertible, with inverse
$
Y\mapsto
\sum_{\ell=0}^\infty ((F_\infty^i)^\top)^\ell
Y(F_\infty^i)^\ell.
$
The derivative is the same for every $K^{-i}$ at $m^i=0$.

Now we perturb around $m^i=0$.
Fix $\widehat K\in\mathcal B_0$. By~\eqref{eq:inf:centered} and~\eqref{eq:inf:dare}, $\Phi^i$ is $C^1$ on the open set where $m^i>-s^i$ and $I+G^i(m^i)P$ is invertible, which contains $(P_\infty^i,\widehat K^{-i},0)$. At this point $\Phi^i$ vanishes and $D_P\Phi^i$ is invertible, so the implicit-function theorem gives neighborhoods $\mathcal V$ of $(\widehat K^{-i},0)$ and $\mathcal W$ of $P_\infty^i$ such that, for every $(K^{-i},m^i)\in\mathcal V$, the equation $\Phi^i(P,K^{-i},m^i)=0$ has exactly one solution $P=P^i(K^{-i};m^i)$ in $\mathcal W$, and this solution is $C^1$ in $(K^{-i},m^i)$.

Next, we glue these local solutions into a single map on $\mathcal B_0$.
The conditions $(s^i+m^i)I-(B^{-i})^\top PB^{-i}\succ0$ and $\rho(F(P,K^{-i},m^i))<1$ hold at $(P_\infty^i,\widehat K^{-i},0)$ by~\eqref{eq:inf:nominal_stabilizing} and are open, so after shrinking $\mathcal V$ they hold at $P=P^i(K^{-i};m^i)$ for every $(K^{-i},m^i)\in\mathcal V$. Since stabilizing solutions of the $\mathcal H_\infty$ Riccati equation are unique~\cite[Ch.~3]{bacsar2008h}, solutions obtained from different $\widehat K$ agree where their neighborhoods overlap. By compactness of $\mathcal B_0$, finitely many of these neighborhoods are needed so that their first coordinate covers $\mathcal{B}_0$ and let $\bar m_i>0$ be the smallest out of all of their second coordinates. This yields the $C^1$ map
$
(K^{-i},m^i)\mapsto P^i(K^{-i};m^i),
\quad
K\in\mathcal B_0, 0\le m^i\le\bar m_i,
$
with $P^i(K^{-i};0)=P_\infty^i$ and
$\Phi^i(P^i(K^{-i};m^i),K^{-i},m^i)=0$.

Shrinking $\bar m_i$ if necessary, we may also assume $\bar Q^i(K^{-i})\succ0$ and $P^i(K^{-i};m^i)\succ0$ for all $K\in\mathcal B_0$ and $m^i\in[0,\bar m_i]$. The first holds since $\bar Q^i(K^{-i})\succeq Q^i-m^i(N-1)\bar k_0^2I$ on $\mathcal B_0$. The second holds since $P^i(K^{-i};0)=P_\infty^i\succ0$ and $P^i$ is uniformly continuous on the compact set $\mathcal B_0\times[0,\bar m_i]$. Writing $F^i(K^{-i};m^i)\coloneqq F(P^i(K^{-i};m^i),K^{-i},m^i)$, we thus have, for all $K\in\mathcal B_0$ and $m^i\in[0,\bar m_i]$,
\begin{subequations}\label{eq:inf:uniform_stabilizing}
\begin{align}
\bar Q^i(K^{-i})&\succ0,\label{eq:inf:us:Q}\\
(s^i+m^i)I-(B^{-i})^\top P^i(K^{-i};m^i)B^{-i}&\succ0,\label{eq:inf:us:margin}\\
\rho\bigl(F^i(K^{-i};m^i)\bigr)&<1.\label{eq:inf:us:schur}
\end{align}
\end{subequations}

Next, we show that $P^i$ is the value of player~$i$'s problem.
By~\eqref{eq:inf:fictitious_dynamics}, \eqref{eq:inf:us:schur} says that the saddle closed loop $A+B^i\mathrm{BR}^i+B^{-i}L^i$ is Schur.
Together with $P^i\succ0$, \eqref{eq:inf:us:Q}, and~\eqref{eq:inf:us:margin}, standard infinite-horizon $\mathcal H_\infty$ verification~\cite[Ch.~3]{bacsar2008h} identifies $P^i$ as the value matrix of~\eqref{eq:inf:centered_problem}, $\mathrm{BR}^i$ as its unique optimal stationary linear controller, and $\widetilde L^i$ as the worst-case disturbance.

Finally, we show that the gains are $C^1$.
The symmetric part of $\mathcal H^i(P^i)$ in~\eqref{eq:inf:saddle} is $\operatorname{diag}\bigl(R^i+(B^i)^\top P^iB^i,\,(s^i+m^i)I-(B^{-i})^\top P^iB^{-i}\bigr)$, which is positive definite by $R^i\succ0$, $P^i\succ0$, and~\eqref{eq:inf:us:margin}. Hence, $\mathcal H^i(P^i)$ is invertible, and since the coefficients of~\eqref{eq:inf:saddle} are $C^1$ in $(K^{-i},m^i)$, so are $\mathrm{BR}^i$ and $L^i$. Their derivatives are continuous on the compact set $\mathcal B_0\times[0,\bar m_i]$, so $C_{ij}<\infty$.

We close this step with a bound used in Step~2.
By~\eqref{eq:complete_square}, \eqref{eq:inf:centered_problem} is the problem with dynamics $x_{t+1}=Ax_t+B^iu_t^i+B^{-i}d_t^i$ and stage payoff $\norm{x_t}^2_{Q^i}+\norm{u_t^i}^2_{R^i}-s^i\norm{d_t^i}^2-m^i\norm{d_t^i-K^{-i}x_t}^2$. For $m^i\ge0$, this payoff is at most that of the decoupled problem ($m^i=0$) at every $(x_t,u_t^i,d_t^i)$. Hence the same ordering holds for the cost of every pair of policies, and thus for the inf--sup values:
\begin{equation}\label{eq:inf:upper}
P^i(K^{-i};m^i)\preceq P_\infty^i
\qquad\text{for all } K\in\mathcal B_0,\ m^i\in[0,\bar m_i].
\end{equation}

\emph{Step 2:}
We now prove part~(c).
Recall from Step~1 that $P^i$, $\mathrm{BR}^i$, and $L^i$ are $C^1$ in $(K^{-i},m^i)$ on $\mathcal B_0\times[0,\bar m_i]$.

To start, we differentiate $P^i$ with respect to $K^j$, $j\neq i$. Let $H^{-i}$ be the perturbation of $K^{-i}$ carrying $H$ in block $j$ and zeros elsewhere, so that $\norm{H^{-i}}=\norm{H}$. By~\eqref{eq:inf:centered},
$D_{K^j}\bar A^i[H]=\frac{m^i}{s^i+m^i}B^{-i}H^{-i}$ and
$D_{K^j}\bar Q^i[H]=-\frac{s^im^i}{s^i+m^i}\big((H^{-i})^\top K^{-i}+(K^{-i})^\top H^{-i}\big)$.
Differentiating $\Phi^i(P^i(K^{-i};m^i),K^{-i},m^i)=0$ in the direction $H$ gives
\begin{equation}\label{eq:inf:implicit}
D_P\Phi^i\bigl[D_{K^j}P^i[H]\bigr]=-\partial_{K^j}\Phi^i[H],
\end{equation}
where $\partial_{K^j}$ differentiates only the explicit dependence on $K^j$, with $P$ and $m^i$ fixed.
To bound $D_{K^j}P^i$, we show that the right-hand side of~\eqref{eq:inf:implicit} is $O(m^i)\norm{H}$ and that $D_P\Phi^i$ has a uniformly bounded inverse.

We start with the right-hand side.
Since $G^i$ does not depend on $K^{-i}$, only $\bar A^i$ and $\bar Q^i$ contribute. Using $P(I+G^iP)^{-1}\bar A^i=PF$, with $F$ as in Step~1, and evaluating at $P=P^i$, we get
\begin{align}\label{eq:inf:partialK}
  \partial_{K^j}\Phi^i[H]
  &=\frac{m^i}{s^i + m^i}\Bigl\{s^i\bigl((H^{-i})^\top K^{-i}+(K^{-i})^\top H^{-i}\bigr)\nonumber\\
  &-(B^{-i}H^{-i})^\top P^iF^i
  -(F^i)^\top P^iB^{-i}H^{-i}\Bigr\}.
\end{align}
We now bound the bracket in~\eqref{eq:inf:partialK}.

The profile $K^{-i}$ appears explicitly in its first term and is bounded by $\sqrt{N-1}\,\bar k_0$ on $\mathcal B_0$; the rest of the dependence on $m^i$ and $K^{-i}$ enters through $P^i$ and $F^i$.
By Step~1 and~\eqref{eq:inf:upper}, $0\prec P^i\preceq P^i_\infty$, so $\norm{P^i}\le\norm{P^i_\infty}$, and it remains to bound $F^i=A+B^i\mathrm{BR}^i+B^{-i}L^i$, which we do through~\eqref{eq:inf:saddle}.
The symmetric part of $\mathcal H^i(P^i)$ is block diagonal, and both of its
blocks are bounded below uniformly:
$R^i+(B^i)^\top P^iB^i\succeq\lambda_{\min}(R^i)I$, while
$P^i\preceq P^i_\infty$ and~\eqref{eq:inf:nominal_stabilizing} give
$(s^i+m^i)I-(B^{-i})^\top P^iB^{-i}\succeq
s^iI-(B^{-i})^\top P^i_\infty B^{-i}\succ0$. Since a matrix whose symmetric
part is at least $\nu I$ has inverse of norm at most $1/\nu$, we get
$\norm{\mathcal H^i(P^i)^{-1}}\le1/\nu^i$ with $\nu^i\coloneqq
\min\{\lambda_{\min}(R^i),\lambda_{\min}(s^iI-(B^{-i})^\top P^i_\infty
B^{-i})\}>0$. As
$P^i\preceq P^i_\infty$, $\norm{K^{-i}}\le\sqrt{N-1}\,\bar k_0$ and $m^i\le\bar m_i$, the
right-hand side $r^i$ of~\eqref{eq:inf:saddle} is bounded, hence so are the
gains $\operatorname{col}(\mathrm{BR}^i,L^i)=\mathcal H^i(P^i)^{-1}r^i$ and, with
them, $F^i=A+B^i\mathrm{BR}^i+B^{-i}L^i$.
Therefore
$\norm{\partial_{K^j}\Phi^i[H]}=O(m^i)\norm{H}$ uniformly on
$\mathcal B_0\times[0,\bar m_i]$.

It remains to bound the inverse of $D_P\Phi^i$.
By~\eqref{eq:inf:us:schur}, $D_P\Phi^i\colon\Delta\mapsto\Delta-(F^i)^\top\Delta F^i$ is invertible at every point of $\mathcal B_0\times[0,\bar m_i]$, so
$
D_{K^j}P^i[H]=-(D_P\Phi^i)^{-1}\,\partial_{K^j}\Phi^i[H].
$
The inverse depends continuously on $(K^{-i},m^i)$ and is therefore bounded on this compact set. Combining the two bounds, $\norm{D_{K^j}P^i}=O(m^i)$ uniformly on $\mathcal B_0\times[0,\bar m_i]$.

Next, we transfer this bound to $\mathrm{BR}^i$. Since $\mathrm{BR}^i$ and $L^i$ are determined jointly by~\eqref{eq:inf:saddle}, let $y^i\coloneqq\operatorname{col}(\mathrm{BR}^i,L^i)$, so that
$\mathcal H^i(P^i)y^i=r^i$.
Differentiating,
$D_{K^j}y^i[H]=\mathcal H^i(P^i)^{-1}\bigl(D_{K^j}r^i[H]-D_{K^j}\mathcal H^i(P^i)[H]\,y^i\bigr)$.
Both $\mathcal H^i$ and $r^i$ depend on $K^{-i}$ only through $P^i$, which enters affinely, and through the term $m^iK^{-i}$ in $r^i$.
By the bound on $D_{K^j}P^i$, both terms in parentheses are therefore $O(m^i)\norm{H}$.
Since $\norm{\mathcal H^i(P^i)^{-1}}\le1/\nu^i$ and $y^i$ is bounded, $\norm{D_{K^j}\mathrm{BR}^i}=O(m^i)$ uniformly on $\mathcal B_0\times[0,\bar m_i]$, i.e., $C_{ij}=O(m^i)$.

For the rate of $K_0^i$, note that 
$K_0^i(m^i)=\mathrm{BR}^i(0;m^i)$ while $K_\infty^i=\mathrm{BR}^i(0;0)$, since at
$m^i=0$ player $i$'s problem is the standalone $\mathcal H_\infty$ problem.
As $\mathrm{BR}^i$ is $C^1$ on the compact set $\mathcal B_0\times[0,\bar m_i]$, it is Lipschitz in $m^i$, so
$K_0^i(m^i)=K_\infty^i+O(m^i)$.

Finally, we verify~\eqref{eq:inf:conditions}. Since $C\ge0$, $\rho(C)\le\max_i\sum_{j\neq i}C_{ij}=O(\max_im^i)$. Moreover, $\norm{K_0^i}+\bar k_0\sum_{j\neq i}C_{ij}\le\norm{K_\infty^i}+O(m^i)$ with $\norm{K_\infty^i}<\bar k_0$. Hence both conditions in~\eqref{eq:inf:conditions} hold for all sufficiently small $m^1,\ldots,m^N$.

\emph{Step 3:}
We now prove part~(b).
First, we bound how each best response depends on the gains. Since $\mathcal B_0$ is convex, for $K,\widetilde K\in\mathcal B_0$, the mean-value inequality and~\eqref{eq:inf:C} give 
\begin{equation}\label{eq:inf:lipschitz}
\norm{\mathrm{BR}^i(K^{-i})-\mathrm{BR}^i(\widetilde K^{-i})}\le\sum_{j\neq i}C_{ij}\norm{K^j-\widetilde K^j}.
\end{equation}

Second, $\mathrm{BR}(\mathcal B_0)\subseteq\mathcal B_0$. For $K\in\mathcal B_0$, \eqref{eq:inf:lipschitz} with $\widetilde K=0$ and~\eqref{eq:inf:conditions} give $\norm{\mathrm{BR}^i(K^{-i})}\le\norm{K_0^i}+\sum_{j\neq i}C_{ij}\norm{K^j}\le\norm{K_0^i}+\bar k_0\sum_{j\neq i}C_{ij}\le\bar k_0$, which gives $\mathrm{BR}(\mathcal B_0)\subseteq\mathcal B_0$.

Third, we show that $\mathrm{BR}$ is a contraction on $\mathcal B_0$ in a weighted norm. Fix $q\in(\rho(C),1)$ and let $w\coloneqq(qI-C)^{-1}\mathbf 1$, where $\mathbf 1$ is the all-ones vector. Since $C\ge0$ and $\rho(C)<q$, the matrix $(qI-C)^{-1}=\sum_{n\ge0}C^n/q^{n+1}$ is entrywise nonnegative, so $w\ge\mathbf 1/q$ and $Cw=qw-\mathbf 1\le qw$. Let $\norm{K}_w\coloneqq\max_i\norm{K^i}/w_i$. Since $\norm{K^j-\widetilde K^j}\le w_j\norm{K-\widetilde K}_w$, \eqref{eq:inf:lipschitz} gives
$
\norm{\mathrm{BR}(K)-\mathrm{BR}(\widetilde K)}_w
\le\max_i\frac{1}{w_i}\sum_{j\neq i}C_{ij}w_j\,\norm{K-\widetilde K}_w
\le q\,\norm{K-\widetilde K}_w.
$

Finally, since $\mathcal B_0$ is closed, Banach's fixed-point theorem yields a unique fixed point $K_\star\in\mathcal B_0$ of $\mathrm{BR}$, with $\norm{\mathrm{BR}^{(n)}(K^{(0)})-K_\star}_w\le q^n\norm{K^{(0)}-K_\star}_w$ for every $K^{(0)}\in\mathcal B_0$. As $q\in(\rho(C),1)$ is arbitrary, the asymptotic root rate is at most $\rho(C)$.
By part~(a), $\mathrm{BR}^i(K^{-i})$ is the unique stationary linear robust best response to $K^{-i}$, so the stationary linear strategically robust dynamic equilibria in $\mathcal B_0$ are exactly the fixed points of $\mathrm{BR}$ in $\mathcal B_0$.
This completes the proof.
\end{proof}

\section*{Acknowledgment}
We used ChatGPT and Claude to assist with proofs (Sections~\ref{sec:open_loop}, \ref{sec:closed_loop},
and~\ref{sec:proofs}) and coding (Section~\ref{sec:experiments}), and to improve grammar, structure, and readability throughout. These tools were not used to formulate the problem, develop research ideas, or draft manuscript sections. We reviewed all AI-assisted work and take responsibility for the final content.
\vspace{-2.0ex}

\section*{References}
\vspace{-3.5ex}
\bibliographystyle{IEEEtran}
\bibliography{references}
\vspace{-9ex}
\begin{IEEEbiography}[{\includegraphics[width=1in,height=1.25in,clip,keepaspectratio]{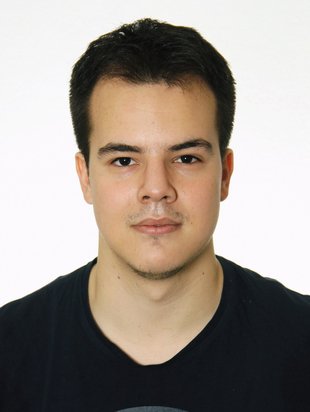}}]{Boris Velasevic} is a Ph.D. student in Computing and Mathematical Sciences at Caltech, advised by Eric Mazumdar. He received the S.B. degrees in mathematics and in computer science and engineering, and the M.Eng.\ degree in computer science and engineering, from MIT\@. His research develops foundations for multi-agent systems, combining control, game theory, and statistics. He received the 2023 MIT Outstanding UROP Student Award and a bronze medal at the 2018 International Mathematical Olympiad.
\end{IEEEbiography}
\vspace{-9ex}
\begin{IEEEbiography}[{\includegraphics[width=1in,height=1.25in,clip,keepaspectratio]{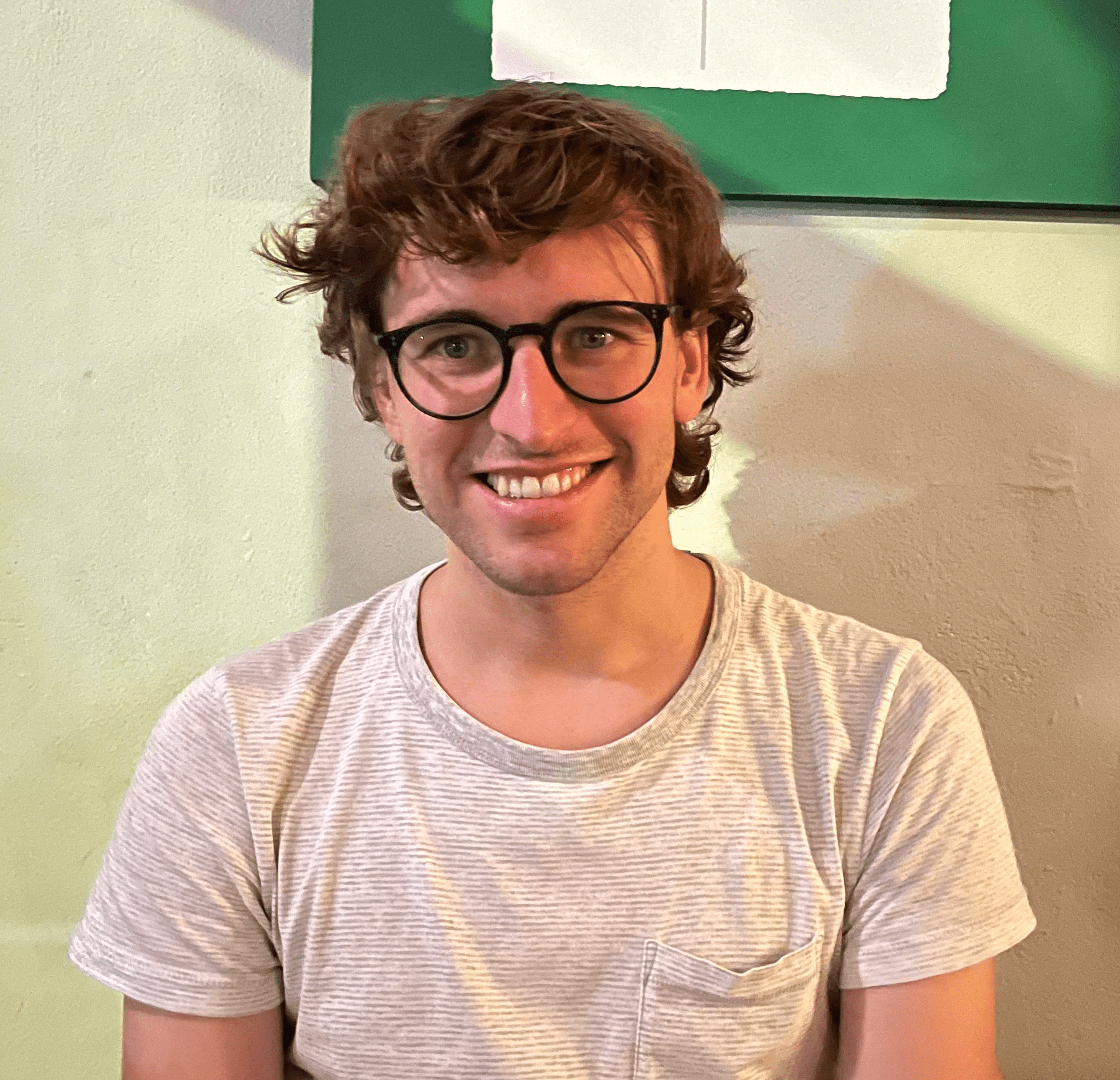}}]{Nicolas Lanzetti} is a PIMCO Postdoctoral Fellow in Data Science at Caltech, working with Adam Wierman, Eric Mazumdar, and Steven Low. He received his Ph.D. from ETH Zurich, advised by Florian D\"{o}rfler (first advisor) and Alessio Figalli (second advisor), and held visiting positions at Stanford, MIT, and EPFL\@. His research develops mathematical and algorithmic foundations for deployable AI in strategic, uncertain, and dynamic environments, using operations research, game theory, optimal transport, machine learning, and control theory.
\end{IEEEbiography}
\vspace{-9ex}
\begin{IEEEbiography}[{\includegraphics[width=1in,height=1.25in,clip,keepaspectratio]{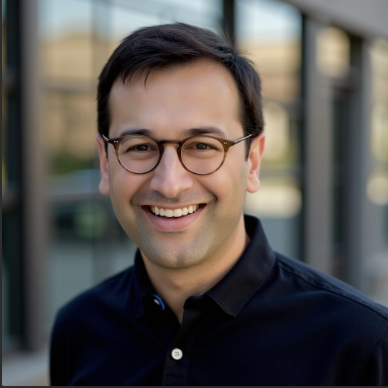}}]{Eric Mazumdar} is an Assistant Professor of Computing and Mathematical Sciences and Economics at Caltech. He received the S.B. degree in electrical engineering and computer science from MIT and the Ph.D. degree in electrical engineering and computer science from UC Berkeley, advised by Michael Jordan and Shankar Sastry. His research lies at the intersection of machine learning, control, and economics, focusing on algorithm design for strategic and dynamic environments. He is a recipient of the NSF CAREER Award, an Okawa Research Award, and a Simons Institute Research Fellowship.
\end{IEEEbiography}

\end{document}

%% file: posting_notice.tex
\AddToShipoutPictureFG*{\AtPageLowerLeft{\put(54,14){\parbox[b]{\dimexpr\paperwidth-108pt\relax}{\fontsize{7}{8}\selectfont This work has been submitted to the IEEE for possible publication. Copyright may be transferred without notice, after which this version may no longer be accessible.}}}}

%% file: figures/meta_figure.tex
\begin{tikzpicture}%
\begin{axis}[%
    width=0.75\columnwidth, height=0.63\columnwidth, scale only axis,%
    xmin=0, xmax=1, ymin=0, ymax=1,%
    xtick={1}, xticklabels={$\infty$}, ytick={1}, yticklabels={$\infty$},%
    tick label style={font=\large}, major tick length=0pt,%
    xlabel={Strategic robustness of P1 $\longrightarrow$},%
    ylabel={Strategic robustness of P2 $\longrightarrow$},%
    label style={font=\normalsize},%
    xlabel style={at={(axis description cs:0.55,-0.04)}, anchor=north},%
    ylabel style={yshift=-4pt},%
    axis line style={line width=0.8pt},%
    axis on top, enlargelimits=false, clip=false,%
    declare function={bA(\x)=0.06*(1-\x)^2; bB(\x)=2*bA(\x)+\x;%
        bnd(\x)=2*bA(\x)/(bB(\x)+sqrt(max(bB(\x)^2-4*bA(\x)^2,0)));},%
]%
    \fill[gray!18] (axis cs:0,0) rectangle (axis cs:1,1);%
    \addplot[draw=none, fill=green!30, domain=0:1, samples=300]%
        {bnd(x)} -- (axis cs:1,1) -- cycle;%
    \addplot[green!55!black, very thick, dash pattern=on 4pt off 2.5pt, domain=0:1, samples=300]%
        {bnd(x)};%
    \fill[black] (axis cs:0,0) circle (2.2pt);%
    \fill[black] (axis cs:1,0) circle (2.2pt);%
    \fill[black] (axis cs:0,1) circle (2.2pt);%
    \node[anchor=north, align=center, font=\small\bfseries, inner sep=0pt]%
        at (axis cs:-0.02,-0.03) {Game-theoretic\\control};%
    \node[anchor=south east, font=\small\bfseries, inner sep=2pt]%
        at (axis cs:0.985,0.024) {P1 robust control};%
    \node[anchor=north west, font=\small\bfseries, inner sep=2pt]%
        at (axis cs:0.018,0.988) {P2 robust control};%
    \draw[-stealth, very thick, black!70, dash pattern=on 4pt off 3pt]%
        (axis cs:0.025,0.025) -- (axis cs:0.9,0.9);%
    \node[rotate=40, anchor=south, font=\small, inner sep=1.5pt,%
          fill=white, fill opacity=0.75, text opacity=1]%
        at (axis cs:0.48,0.5) {strategically robust equilibria};%
    \node[align=center, font=\normalsize, text=green!25!black]%
        at (axis cs:0.76,0.40) {unique \&\\ tractable\\ equilibria};%
\end{axis}%
\end{tikzpicture}%

%% file: figures/free_lunch_plot.tex
\definecolor{cbBlue}{HTML}{0072B2}
\definecolor{cbOrange}{HTML}{E69F00}
\definecolor{cbGreen}{HTML}{009E73}
\definecolor{cbVermillion}{HTML}{D55E00}
\def\data{data/cost_vs_M_2p_asymmetric.csv}
\def\datasym{data/cost_vs_M_2p_symmetric.csv}
\begin{tikzpicture}
\begin{axis}[
    width=7.7cm,
    height=4.2cm,
    xlabel={$1/m$},
    ylabel={Social cost},
    xmin=0.01,
    xmax=5.0,
    ymin=0.12,
    ymax=0.55,
    xtick={0,1,2,3,4,5},
    ytick={0.15, 0.2, 0.3, 0.4, 0.5},
    clip=true,
    grid=major,
    grid style={gray!45},
    legend style={at={(0.97,0.97)}, anchor=north east, font=\footnotesize, draw=gray!50, fill=white, fill opacity=0.92, text opacity=1},
    legend cell align={left},
    legend image post style={xscale=0.85},
    tick label style={font=\footnotesize},
    label style={font=\small},
    after end axis/.code={
        \draw[->, line width=0.9pt]
            (axis description cs:0.32,-0.18) -- (axis description cs:0.04,-0.18);
        \node[anchor=north, font=\footnotesize]
            at (axis description cs:0.18,-0.18) {less robust};
        \draw[->, line width=0.9pt]
            (axis description cs:0.68,-0.18) -- (axis description cs:0.96,-0.18);
        \node[anchor=north, font=\footnotesize]
            at (axis description cs:0.82,-0.18) {more robust};
    }
]
\addplot[color=cbBlue, line width=1.5pt]
    table [col sep=comma, x={1/M}, y=TotalCost] {\data};
\addlegendentry{Player 1 robust}
\addplot[color=cbOrange, line width=1.5pt]
    table [col sep=comma, x={1/M}, y=TotalCost] {\datasym};
\addlegendentry{Both players robust}
\addplot[black, dashed, line width=1.6pt]
    coordinates {(\pgfkeysvalueof{/pgfplots/xmin},0.153846) (\pgfkeysvalueof{/pgfplots/xmax},0.153846)};
\addlegendentry{Social optimum}
\end{axis}
\end{tikzpicture}

%% file: figures/free_lunch_plot_network.tex
\definecolor{cbBlue}{HTML}{0072B2}
\definecolor{cbOrange}{HTML}{E69F00}

\def\data{data/cost_vs_M_cent_asymmetric.csv}
\def\datasym{data/cost_vs_M_cent_symmetric.csv}

\begin{tikzpicture}
\begin{axis}[
    width=7.7cm,
    height=4.2cm,
    xlabel={$1/m$},
    ylabel={Social cost},
    xmin=0,
    xmax=40,
    ymin=3.8,
    ymax=14,
    clip mode=individual, %
    grid=major,
    grid style={gray!45},
    legend style={at={(0.5,0.97)}, anchor=north, font=\footnotesize, draw=gray!50, fill=white, fill opacity=0.92, text opacity=1},
    legend cell align={left},
    legend image post style={xscale=0.85},
    tick label style={font=\footnotesize},
    label style={font=\small},
]

\addplot[color=cbBlue, line width=1.5pt]
    table [col sep=comma, x={1/M}, y=TotalCost] {\data};
\addlegendentry{Center robust}

\addplot[color=cbOrange, line width=1.5pt]
    table [col sep=comma, x={1/M}, y=TotalCost] {\datasym};
\addlegendentry{All players robust}

\addplot[black, dashed, line width=1.6pt]
    coordinates {(\pgfkeysvalueof{/pgfplots/xmin},4.011799) (\pgfkeysvalueof{/pgfplots/xmax},4.011799)};
\addlegendentry{Social optimum}

\draw[->, line width=0.9pt]
    (axis description cs:0.32,-0.18) -- (axis description cs:0.04,-0.18);
\node[anchor=north, font=\footnotesize]
    at (axis description cs:0.18,-0.18) {less robust};

\draw[->, line width=0.9pt]
    (axis description cs:0.68,-0.18) -- (axis description cs:0.96,-0.18);
\node[anchor=north, font=\footnotesize]
    at (axis description cs:0.82,-0.18) {more robust};

\end{axis}
\end{tikzpicture}

%% file: figures/platoon_illustration.tex
\begin{tikzpicture}[
    scale=0.85,
    every node/.style={font=\footnotesize},
    car/.pic={
        \draw[line width=0.8pt, fill=gray!15, rounded corners=2pt]
            (-0.6,0.12) -- (-0.6,0.44) -- (-0.32,0.44) -- (-0.18,0.70)
            -- (0.24,0.70) -- (0.40,0.44) -- (0.6,0.44) -- (0.6,0.12) -- cycle;
        \draw[line width=0.8pt, fill=white] (-0.34,0.12) circle (0.11) (0.34,0.12) circle (0.11);
    },
    disp/.style={-stealth, line width=0.8pt},
    gap/.style={stealth-stealth, line width=0.6pt},
]
\draw[gray!60, line width=0.6pt] (-0.9,0) -- (7.8,0);

\foreach \i/\x in {4/0, 3/2.3, 2/4.6, 1/6.9} {
    \pic at (\x,0) {car};
    \node[font=\scriptsize] at (\x,0.29) {$\i$};
    \draw[disp] (\x-0.25,1.0) -- (\x+0.35,1.0);
    \node[anchor=south] at (\x+0.05,1.02) {$u_t^{\i}$};
}

\foreach \i/\x in {4/0, 3/2.3, 2/4.6} {
    \draw[gap] (\x+0.6,-0.22) -- (\x+1.7,-0.22);
    \node[anchor=north] at (\x+1.15,-0.25) {$g_t^{\i}$};
}
\node[anchor=north, font=\scriptsize] at (6.9,-0.12) {leader};
\end{tikzpicture}

%% file: figures/platoon_trajectories.tex
\definecolor{cbBlue}{HTML}{0072B2}%
\definecolor{cbVermillion}{HTML}{D55E00}%
\definecolor{cbGreen}{HTML}{009E73}%
\def\agA{data/platoon_gaps_m160_gap5.csv}%
\def\agB{data/platoon_gaps_m1000_gap5.csv}%
\def\bgA{data/platoon_gaps_m160_gap20.csv}%
\def\bgB{data/platoon_gaps_m1000_gap20.csv}%
\def\gdes{20}%
\begin{tikzpicture}
\begin{groupplot}[
  group style={group size=2 by 1, horizontal sep=0.45cm, y descriptions at=edge left},
  scale only axis, height=2.9cm, set layers,
  xlabel={Time step $t$}, ylabel={Gap $g_t^i$},
  xmin=0, xmax=20, ymin=0, ymax=25, grid=major, grid style={gray!40},
  tick label style={font=\footnotesize}, label style={font=\small},
  title style={font=\small}, no markers,
  legend cell align={left},
]
\nextgroupplot[width=4.55cm, title={(a) $g_0^i=5$}]
  \pgfonlayer{axis background}
    \fill[gray!14] (axis cs:0,0) rectangle (axis cs:20,\gdes);
  \endpgfonlayer
  \draw[black, densely dotted] (axis cs:0,\gdes) -- (axis cs:20,\gdes);
  \node[anchor=south east, font=\scriptsize, text=black!60] at (axis cs:19.6,0.4) {safety penalty active};
  \addplot[cbBlue,       line width=1.1pt] table[col sep=comma,x=t,y=g2]{\agA};
  \addplot[cbVermillion, line width=1.1pt] table[col sep=comma,x=t,y=g3]{\agA};
  \addplot[cbGreen,      line width=1.1pt] table[col sep=comma,x=t,y=g4]{\agA};
  \addplot[cbBlue,       dashed, line width=1.1pt] table[col sep=comma,x=t,y=g2]{\agB};
  \addplot[cbVermillion, dashed, line width=1.1pt] table[col sep=comma,x=t,y=g3]{\agB};
  \addplot[cbGreen,      dashed, line width=1.1pt] table[col sep=comma,x=t,y=g4]{\agB};
\nextgroupplot[width=2.4cm, title={(b) $g_0^i=20$}, xtick={0,10,20}, xlabel={},
  legend style={at={(0.5,0.04)}, anchor=south, font=\scriptsize, draw=gray!50,
                fill=white, fill opacity=0.9, text opacity=1, row sep=-1.5pt},
  legend image post style={xscale=0.8},
]
  \pgfonlayer{axis background}
    \fill[gray!14] (axis cs:0,0) rectangle (axis cs:20,\gdes);
  \endpgfonlayer
  \draw[black, densely dotted] (axis cs:0,\gdes) -- (axis cs:20,\gdes);
  \addplot[cbGreen,      dashed, line width=1.1pt, forget plot] table[col sep=comma,x=t,y=g4]{\bgB};
  \addplot[cbVermillion, dashed, line width=1.1pt, forget plot] table[col sep=comma,x=t,y=g3]{\bgB};
  \addplot[cbBlue,       dashed, line width=1.1pt, forget plot] table[col sep=comma,x=t,y=g2]{\bgB};
  \addplot[cbGreen,      line width=1.1pt, forget plot] table[col sep=comma,x=t,y=g4]{\bgA};
  \addplot[cbVermillion, line width=1.1pt, forget plot] table[col sep=comma,x=t,y=g3]{\bgA};
  \addplot[cbBlue,       line width=1.1pt, forget plot] table[col sep=comma,x=t,y=g2]{\bgA};
  \node[anchor=south, font=\scriptsize, align=center] at (axis cs:10,20.6) {all gaps, both $m$};
  \addlegendimage{cbBlue, line width=1.1pt}        \addlegendentry{vehicle 2}
  \addlegendimage{cbVermillion, line width=1.1pt}  \addlegendentry{vehicle 3}
  \addlegendimage{cbGreen, line width=1.1pt}       \addlegendentry{vehicle 4}
  \addlegendimage{black, line width=1.1pt}         \addlegendentry{$m=160$}
  \addlegendimage{black, dashed, line width=1.1pt} \addlegendentry{$m=1000$}
\end{groupplot}
\end{tikzpicture}%

%% file: figures/tullock_plot.tex
\definecolor{cbBlue}{HTML}{0072B2}
\definecolor{cbGreen}{HTML}{009E73}
\definecolor{cbOrange}{HTML}{E69F00}
\definecolor{cbVermillion}{HTML}{D55E00}
\def\dataA{data/tullock_contest_m0.05_P1_controls.csv}
\def\dataB{data/tullock_contest_m0.1_P1_controls.csv}
\def\dataC{data/tullock_contest_m0.2_P1_controls.csv}
\def\dataD{data/tullock_contest_m0.4_P1_controls.csv}
\begin{tikzpicture}
\begin{axis}[
    width=7.7cm,
    height=4.2cm,
    xlabel={Time step $t$},
    ylabel={Player acquisition $u_t^i$},
    xmin=-0.3, xmax=9.3,
    ymin=-0.18, ymax=1.2,
    xtick={0,1,2,3,4,5,6,7,8,9},
    ytick={0,0.5,1},
    grid=major,
    grid style={gray!45},
    legend style={at={(0.03,0.97)}, anchor=north west, font=\footnotesize,
                  draw=gray!50, fill=white, fill opacity=0.92, text opacity=1},
    legend cell align={left},
    legend columns=2,
    tick label style={font=\footnotesize},
    label style={font=\small},
    every axis plot/.append style={const plot, line width=1.3pt},
    no markers,
]

\addplot[color=cbBlue]       table[col sep=comma, x=t, y=u_P1] {\dataA};
\addlegendentry{$m=0.05$}
\addplot[color=cbGreen]      table[col sep=comma, x=t, y=u_P1] {\dataB};
\addlegendentry{$m=0.1$}
\addplot[color=cbOrange]     table[col sep=comma, x=t, y=u_P1] {\dataC};
\addlegendentry{$m=0.2$}
\addplot[color=cbVermillion] table[col sep=comma, x=t, y=u_P1] {\dataD};
\addlegendentry{$m=0.4$}
\end{axis}
\end{tikzpicture}

%% file: figures/cournot_plot.tex
\definecolor{cbBlue}{HTML}{0072B2}
\definecolor{cbGreen}{HTML}{009E73}
\definecolor{cbOrange}{HTML}{E69F00}
\definecolor{cbVermillion}{HTML}{D55E00}
\def\dataA{data/resource_game_m1.0_P1_controls.csv}
\def\dataB{data/resource_game_m2.0_P1_controls.csv}
\def\dataC{data/resource_game_m3.0_P1_controls.csv}
\def\dataD{data/resource_game_m5.0_P1_controls.csv}
\begin{tikzpicture}
\begin{axis}[
    width=7.7cm,
    height=4.2cm,
    xlabel={Time step $t$},
    ylabel={Player acquisition $u_t^i$},
    xmin=-0.3, xmax=9.3,
    ymin=-0.06, ymax=1.14,
    xtick={0,1,2,3,4,5,6,7,8,9},
    ytick={0,0.5,1},
    grid=major,
    grid style={gray!45},
    legend style={at={(0.03,0.97)}, anchor=north west, font=\footnotesize,
                  draw=gray!50, fill=white, fill opacity=0.92, text opacity=1},
    legend cell align={left},
    legend columns=2,
    tick label style={font=\footnotesize},
    label style={font=\small},
    every axis plot/.append style={const plot, line width=1.3pt},
    no markers,
]
\addplot[color=cbBlue]       table[col sep=comma, x=t, y=u_P1] {\dataA};
\addlegendentry{$m=1.0$}
\addplot[color=cbGreen]      table[col sep=comma, x=t, y=u_P1] {\dataB};
\addlegendentry{$m=2.0$}
\addplot[color=cbOrange]     table[col sep=comma, x=t, y=u_P1] {\dataC};
\addlegendentry{$m=3.0$}
\addplot[color=cbVermillion] table[col sep=comma, x=t, y=u_P1] {\dataD};
\addlegendentry{$m=5.0$}
\end{axis}
\end{tikzpicture}